\magnification=\magstep1%
\input amstex
\input pictex
\input xy
\xyoption{all}
\documentstyle{amsppt}
\parskip=4pt
\hfuzz=12pt
\NoBlackBoxes
\nologo
\def\la#1{\hbox to #1pc{\leftarrowfill}}
\def\ra#1{\hbox to #1pc{\rightarrowfill}}

\nopagenumbers

\define\R{{\Bbb R}}
\define\Z{{\Bbb Z}}
\define\GT{{\Cal G  \Cal T}}
\define\GCCT{{\Cal G \frak C \Cal C \Cal T}}

\topmatter
\title The algebraic construction of the Novikov complex of a circle-valued 
Morse function
\endtitle
\rightheadtext{}
\leftheadtext{}
\abstract{The Novikov complex of a circle-valued Morse function $f:M
\to S^1$ is constructed algebraically from the Morse-Smale complex of
the restriction of the real-valued Morse function
$\overline{f}:\overline{M} \to \R$
to a fundamental domain of the pullback infinite cyclic
cover $\overline{M}=f^*\R$ of $M$.}
\endabstract
\author Andrew Ranicki\endauthor
\address Department of Mathematics and Statistics \newline
\indent University of Edinburgh\newline
\indent Edinburgh EH9 3JZ\newline 
\indent Scotland, UK\newline
\vbox{\vskip1mm}
\endaddress

\email  aar\@maths.ed.ac.uk \endemail
\endtopmatter
\nopagenumbers

\subheading{Introduction}
\footnote[~]{To appear in Mathematische Annalen}

The relationship between real-valued Morse functions and chain
complexes is well understood.  The Morse-Smale complex of a Morse
function $f:M \to \R$ on a compact $m$-dimensional manifold $M$ is
defined  using a choice of gradient-like vector field $v$ 
satisfying the transversality condition, to be a based f.\,g.  
free $\Z[\pi_1(M)]$-module chain complex $C^{MS}(M,f,v)$ with
$$\hbox{\rm rank}_{\Z[\pi_1(M)]}C^{MS}_i(M,f,v)~=~c_i(f)$$
the number of critical points of $f$ with index $i$, and the differentials
defined by counting the downward $v$-gradient flow lines 
in the universal cover $\widetilde{M}$ of $M$.  
The pair $(f,v)$ determines a handlebody decomposition on $M$
with one $i$-handle for each critical point of index $i$
$$M~=~\bigcup\limits^m_{i=0}\bigcup_{c_i(f)}D^i\times D^{m-i}~,$$
and the Morse-Smale complex is the cellular chain complex of $\widetilde{M}$
$$C^{MS}(M,f,v)~=~C(\widetilde{M})~.$$
\indent The relationship between circle-valued Morse functions $f:M \to S^1$
and chain complexes is more complicated, and not so well understood.
A lift of $f$ to the pullback infinite cyclic cover $\overline{M}=f^*\R$
is a real-valued Morse function $\overline{f}:\overline{M} \to \R$
on a non-compact manifold, so traditional Morse theory does not
apply directly. The methods developed (by Novikov, Farber, Pajitnov, the author and others)
to count the critical points of $f$ use the structure of the fundamental
group ring $\Z[\pi_1(M)]$ as a Laurent polynomial extension of
$\Z[\pi_1(\overline{M})]$, 
as well as a completion and a localization of $\Z[\pi_1(M)]$.
For any map $f:M \to S^1$ with $M$ and
$\overline{M}$ connected the ring $\Z[\pi_1(M)]$
is the $\alpha $-twisted Laurent polynomial extension 
$$\Z[\pi_1(M)]~=~\Z[\pi_1(\overline{M})]_{\alpha }[z,z^{-1}]~
=~\Z[\pi_1(\overline{M})]_{\alpha }[z][z^{-1}]$$
with $\alpha :\pi_1(\overline{M}) \to \pi_1(\overline{M})$ the monodromy
automorphism and
$$az~=~z\alpha (a)~~(a \in\Z[\pi_1(\overline{M})])~.$$
The completion
$$\widehat{\Z[\pi_1(M)]}~=~\Z[\pi_1(\overline{M})]_{\alpha }((z))~
=~\Z[\pi_1(\overline{M})]_{\alpha }[[z]][z^{-1}]$$
is called the {\it Novikov ring} of $\Z[\pi_1(M)]$. 
Let $\Sigma$ be the set of square matrices in 
$\Z[\pi_1(\overline{M})]_{\alpha }[z] \subset \Z[\pi_1(M)]$ which become 
invertible in $\Z[\pi_1(\overline{M})]$ under the augmentation $z \mapsto 0$.
The noncommutative localization $\Sigma^{-1}\Z[\pi_1(M)]$ of $\Z[\pi_1(M)]$
(in the sense of Cohn \cite{C}) is a ring with a morphism
$\Z[\pi_1(M)] \to \Sigma^{-1}\Z[\pi_1(M)]$
such that any ring morphism $\Z[\pi_1(M)] \to R$ which sends $\Sigma$
to invertible matrices in $R$ has a unique factorization
$\Z[\pi_1(M)] \to \Sigma^{-1}\Z[\pi_1(M)] \to R$. The inclusion 
$\Z[\pi_1(M)] \to \widehat{\Z[\pi_1(M)]}$ sends $\Sigma$ to
invertible matrices in $\widehat{\Z[\pi_1(M)]}$, so there is a 
natural ring morphism $\Sigma^{-1}\Z[\pi_1(M)]\to\widehat{\Z[\pi_1(M)]}$.

A vector field on  a manifold $M$ is a section of the tangent bundle of $M$
$$v~:~M \to \tau_M~.$$ 
The gradient of a Morse function $f:M \to \R$ is a section
of the cotangent bundle 
$$\nabla f~=~(\partial f /\partial x_i) ~:~M \to \tau_M^*~=~\bigcup_{x \in M} 
\hbox{\rm Hom}_{\R}(\tau_M(x),\R)$$
with zeros the critical points of $f$. 
A vector field $v:M \to \tau_M$ is 'gradient-like' for $f$ if
there exists a Riemannian metric $\beta :\tau_M \cong \tau^*_M$ 
on $M$ such that 
$$\beta  \circ v~=~\nabla f ~:~M \to \tau_M^*~.$$
A $v$-gradient flow line $\gamma: \R \to M$ satisfies
$$\gamma'(t)~=~-v(\gamma(t))  \in \tau_M(\gamma(t))~.$$
The limits
$$\lim_{t \to -\infty} \gamma(t)~=~p~~,~~\lim_{t \to \infty} \gamma(t)~=~q \in M$$
are critical points of $f$. For every non-critical point $x \in M$ there is a
$v$-gradient flow line $\gamma_x:\R \to M$ (which is unique up to scaling)
such that $\gamma_x(0)=x \in M$. The unstable and stable manifolds of 
a critical point $p \in M$ of index $i$ are the open manifolds
$$\aligned
&W^u(p,v)~=~\{x\in M\,\vert\,\lim_{t\rightarrow -\infty}\gamma_x(t)=p\}~,\cr
&W^s(p,v) ~=~\{x\in M\,\vert\,\lim_{t\rightarrow\infty}\gamma_x(t)=p\}
\endaligned$$
which are diffeomorphic to $\R^i$, $\R^{m-i}$ respectively. 
Let $\GT(f)$ denote the space of all gradient-like vector fields $v$ for $f$
which satisfy the Morse-Smale transversality condition that 
for any distinct critical points $p,q \in M$ the submanifolds
$W^u(p,v),W^s(q,v) \subset M$ intersect transversely. 
For $v \in \GT(f)$ and critical points $p,q \in M$ of index $i,i-1$ 
there is only a finite number of $v$-gradient flow lines $\gamma : \R \to M$
which start at $p$ and terminate at $q$, and choosing orientations
there is obtained an  algebraic number $n(p,q)\in \Z$.

A circle-valued Morse function $f:M \to S^1$ lifts to a real-valued
Morse function $\overline{f}:\overline{M} \to \R$ on the infinite
cyclic cover $\overline{M}=f^*\R$. Let $\GT(f)$ be the space of all
vector fields $v$ on $M$ which lift to a gradient-like vector
field $\overline{v}$ on $\overline{M}$ satisfying the transversality
condition, so that $\overline{v} \in \GT(\overline{f})$.

Novikov \cite{N} and Pajitnov \cite{P1} used the completion
$\widehat{\Z[\pi_1(M)]}$ to construct geometrically for
any circle-valued Morse function $f:M \to S^1$ and 
$v \in \GT(f)$ a based f.\,g. free $\widehat{\Z[\pi_1(M)]}$-module 
chain complex $C^{Nov}(M,f,v)$ such that 
$$\hbox{\rm rank}_{\widehat{\Z[\pi_1(M)]}}C^{Nov}_i(M,f,v)~=~c_i(f)$$
with $c_i(f)$ the number of critical points of $f$ with index $i$.  
As in the real-valued case the differentials are defined by counting the
$\widetilde{v}$-gradient flow lines in the universal
cover $\widetilde{M}$.  Moreover, there is a chain equivalence
$$C^{Nov}(M,f,v)~\simeq~C(\widetilde{M};\widehat{\Z[\pi_1(M)]})$$
with 
$$C(\widetilde{M};\widehat{\Z[\pi_1(M)]})~=~
\widehat{\Z[\pi_1(M)]}\otimes_{\Z[\pi_1(M)]}C(\widetilde{M})$$
the $\widehat{\Z[\pi_1(M)]}$-coefficient cellular chain complex of 
$\widetilde{M}$, for any $CW$ structure on $M$.

Pajitnov \cite{P2,P3,P4} constructed for any Morse function
$f:M \to S^1$ a $C^0$-dense subspace  $\GCCT(f)\subset \GT(f)$, such
that for $v \in \GCCT(f)$ the coefficients in the corresponding Novikov 
complex $C^{Nov}(M,f,v)$ are rational, in the sense that 
$$C^{Nov}(M,f,v)~=~
\widehat{\Z[\pi_1(M)]}\otimes_{\Sigma^{-1}\Z[\pi_1(M)]}C^{Paj}(M,f,v)$$
for a based f.\,g.  free  $\Sigma^{-1}\Z[\pi_1(M)]$-module chain 
complex $C^{Paj}(M,f,v)$, with a chain equivalence
$$C^{Paj}(M,f,v)~\simeq~C(\widetilde{M};\Sigma^{-1}\Z[\pi_1(M)])~.$$
(Strictly speaking, $C^{Paj}(M,f,v)$ was only defined for abelian $\pi_1(M)$,
when the natural ring morphism
$\Sigma^{-1}\Z[\pi_1(M)] \to \widehat{\Z[\pi_1(M)]}$ is injective,
but this was for algebraic convenience rather than out of geometric necessity.) 

For a Morse function $f:M \to S^1$ which is transverse regular at
$0 \in S^1$ the lift $\overline{f}:\overline{M} \to \R$ 
is transverse regular at $\Z \subset \R$, and
the restriction to a fundamental domain is a Morse function 
$$f_N~=~\overline{f}\vert ~:~
(M_N;N,N_1)~=~\overline{f}^{\,-1}(I;\{0\},\{1\}) \to (I;\{0\},\{1\})$$
with $c_i(f_N)=c_i(f)$ critical points of index $i$, and $N_1$ a copy of $N$.
Every $v \in \GT(f)$ lifts to $\overline{v} \in \GT(\overline{f})$,
and $v_N=\overline{v}\vert \in \GT(f_N)$ determines a handlebody decomposition
$$M_N~=~N \times I \cup \bigcup\limits^m_{i=0}\bigcup\limits_{c_i(f)}D^i \times D^{m-i}$$
with one $i$-handle for each index $i$ critical point of $f$.
Given a $CW$ structure on $N$ with $c_i(N)$ $i$-cells use the handlebody
structure on $M_N$ to define  a $CW$ structure on $M_N$
with  $c_i(N)+c_i(f)$ $i$-cells.  (In practice, the $CW$ structure will
be the one determined by a Morse function $N \to \R$).
The inclusion $g:N \to M_N$ induces an inclusion of based
f.g. free $\Z[\pi_1(\overline{M})]$-module chain complexes
$$g~:~C(\widetilde{N}) \to C(\widetilde{M}_N)$$
with $\widetilde{N},\widetilde{M}_N$ the covers of $N,M_N$ induced
from the universal cover $\widetilde{M}$ of $\overline{M}$.
Write the inclusion of $N_1$ in $M_N$ as $h:N \to M_N$,
and note that in general $N_1=h(N)$ is not a $CW$ subcomplex of $M_N$.
A {\it chain approximation} for $h$ is a
$\Z[\pi_1(\overline{M})]$-module chain map
$$h~:~C(\widetilde{N}_1)~=~\alpha C(\widetilde{N}) \to C(\widetilde{M}_N)$$
in the chain homotopy class induced by the map $h:N \to M_N$
as given by the cellular approximation theorem.
For any choice of chain approximation $h$ Farber and Ranicki \cite{FR} 
defined algebraically a based f.\,g. free $\Sigma^{-1}\Z[\pi_1(M)]$-module 
chain complex $C^{FR}(M,f,v,h)$ such that
$$\hbox{\rm rank}_{\Sigma^{-1}\Z[\pi_1(M)]}C^{FR}_i(M,f,v,h)~=~c_i(f)$$
as a deformation of the
$\Sigma^{-1}\Z[\pi_1(M)]$-coefficient Morse-Smale complex 
$$C^{MS}(M_N,f_N,v_N;\Sigma^{-1}\Z[\pi_1(M)])~=~\Sigma^{-1}\Z[\pi_1(M)]
\otimes_{\Z[\pi_1(\overline{M})]}C^{MS}(M_N,f_N,v_N)$$
with a chain equivalence
$$C^{FR}(M,f,v,h)~\simeq~\Sigma^{-1}\Z[\pi_1(M)]\otimes_{\Z[\pi_1(M)]}C(\widetilde{M})~.$$ 
There is also a $\widehat{\Z[\pi_1(M)]}$-coefficient version
$$C^{FR}(M,f,v,h;\widehat{\Z[\pi_1(M)]})~=~
\widehat{\Z[\pi_1(M)]}\otimes_{\Sigma^{-1}\Z[\pi_1(M)]}C^{FR}(M,f,v,h)~.$$

\proclaim{\bf Cokernel Theorem 6.6} The chain complex of \cite{FR} is
isomorphic to the cokernel
$$C^{FR}(M,f,v,h)~\cong~\hbox{\rm coker}(\phi)$$
of the morphism of based f.g. free $\Sigma^{-1}\Z[\pi_1(M)]$-module
chain complexes 
$$\phi~=~g-zh~:~\Sigma^{-1}C(\widetilde{N})_{\alpha}[z,z^{-1}] \to 
\Sigma^{-1}C(\widetilde{M}_N)_{\alpha}[z,z^{-1}]$$
which is a split injection in each degree.\hfill\qed
\endproclaim

The expression $C^{FR}(M,f,v,h)$ as a cokernel makes it possible
to prove invariance results such as :

\proclaim{\bf Invariance Theorem 6.7} 
Let $f:M \to S^1$ be a Morse function, and let $v \in \GT(f)$.
For any regular values $0,0' \in S^1$, 
$CW$ structures on $N=f^{-1}(0)$, $N'=f^{-1}(0')$ and chain approximations
$$h~:~\alpha C(\widetilde{N}) \to C(\widetilde{M}_N)~,~
h'~:~\alpha C(\widetilde{N}') \to C(\widetilde{M}_{N'})$$
there is defined a simple isomorphism of 
based f.\,g. free $\Sigma^{-1}\Z[\pi_1(M)]$-module chain complexes
$$C^{FR}(M,f,v,h)~\cong~C^{FR}(M,f,v,h')~.\eqno{\qed}$$
\endproclaim

Here, simple means that the torsion of the isomorphism is in the image of
$\{\pm \pi_1(M)\} \subseteq K_1(\Sigma^{-1}\Z[\pi_1(M)])$.

Given $f:M \to S^1$, $v \in \GT(f)$ and a $CW$ structure on 
$N=f^{-1}(0) \subset M$ we shall say that a chain approximation
$$h^{gra}~:~\alpha C(\widetilde{N}) \to C(\widetilde{M}_N)$$
is  {\it gradient-like} if it counts the $v$-gradient flow lines
in the universal cover $\widetilde{M}$. (The precise definition is given
in 6.8).

\noindent{\bf Isomorphism Theorem 6.10} {\it
For $v \in \GT(f)$ with a gradient-like chain approximation $h^{gra}$
there are basis-preserving isomorphisms
$$\aligned
&C^{Nov}(M,f,v)~\cong~C^{FR}(M,f,v,h^{gra};\widehat{\Z[\pi_1(M)]})~,\cr
&C^{Paj}(M,f,v)~\cong~C^{FR}(M,f,v,h^{gra})~.\hfill\qed
\endaligned$$}

Pajitnov \cite{P3} (Theorem 7.2) showed that every $v \in \GT(f)$
is $C^0$-close to $v^{gra} \in \GCCT(f)$ with a gradient-like chain approximation 
$h^{gra}:C(\widetilde{N}) \to C(\widetilde{M}_N)$.
A similar construction was obtained by Hutchings and Lee \cite{HL}.
In fact, Cornea and Ranicki \cite{CR} prove that every $v \in \GT(f)$
admits a gradient-like chain approximation $h^{gra}$. 

The plan of the paper is as follows.  \S1 is purely algebraic in
nature, concerning the chain homotopy properties of the algebraic
mapping cones and cokernels of chain maps.  The glueing
properties of the Morse-Smale complex $C^{MS}(M,f,v)$ 
are described in \S2 for finite unions, and in \S3 for infinite unions.
\S4 gives a brief account of the Cohn localization.
\S5 deals with the cokernel and infinite union
construction of chain complexes over a twisted polynomial extension
$A_{\alpha }[z,z^{-1}]$ for any ring $A$ with automorphism $\alpha :A \to
A$, the localization $\Sigma^{-1}A_{\alpha }[z,z^{-1}]$ and the Novikov
ring $A_{\alpha }((z))=A_{\alpha }[[z]][z^{-1}]$. The Cokernel,
Invariance and Isomorphism Theorems are proved in \S6.

The remainder of the Introduction is an outline of the proof of the
Cokernel, Invariance and Isomorphism Theorems in the special case when
$f_*:\pi_1(M) \to \pi_1(S^1)$ is an isomorphism, so that
$$\aligned
&\pi_1(\overline{M})~=~\{1\}~,~\Z[\pi_1(M)]~=~\Z[z,z^{-1}]~,\cr
&\Sigma^{-1}\Z[\pi_1(M)]~=~(1+z\Z[z])^{-1}\Z[z,z^{-1}]~,\cr
&\widehat{\Z[\pi_1(M)]}~=~\Z((z))~.
\endaligned$$
\indent 
The Novikov complex of a Morse function $f:M \to S^1$ with 
respect to $v \in \GT(f)$ is the based f.\,g. free $\Z((z))$-module 
chain complex 
$$C~=~C^{Nov}(M,f,v)$$ 
with
$$d_C~:~C_i~=~\Z((z))^{c_i(f)} \to C_{i-1}~=~\Z((z))^{c_{i-1}(f)}~;~
\overline{p} \mapsto \sum\limits_{\overline{q}}
 n(\overline{p},\overline{q})\overline{q}$$
where $n(\overline{p},z^j\overline{q})$ is the algebraic number of
$\overline{v}$-gradient flow lines in $\overline{M}$ from a critical point
$\overline{p} \in \overline{M}$ of index $i$ to a critical point
$\overline{q} \in \overline{M}$ of index $i-1$, using the
transversality property of $v$ to ensure that these numbers are finite.
{\it Here $z:\overline{M} \to \overline{M}$ is the generating covering 
translation parallel to the $v$-gradient flow of $f$.} 
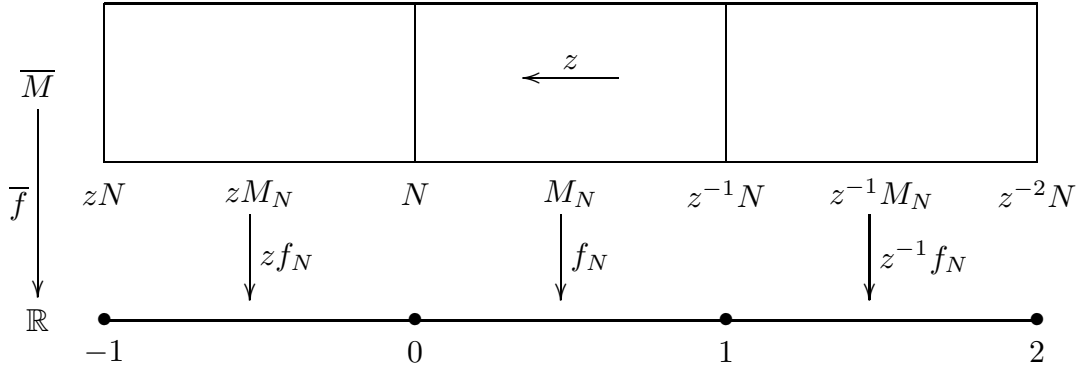
$$\beginpicture
\setcoordinatesystem units <7pt,10pt>  
\put {$\overline{M}$} at -24 -2.5
\put {$\xymatrix@R+35pt {\ar[d]_-{\displaystyle \overline{f}} & \\ &}$} 
at -22.5 -6.3
\put {$\R$} at -24 -10
\putrule from -21 0 to 21 0
\putrule from -7 0 to -7 -5
\putrule from -21 0 to -21 -5
\putrule from 21 0 to 21 -5
\putrule from 7 0 to 7 -5
\putrule from -21 -5 to 21 -5
\put {$zM_N$} at -14 -6
\put {$M_N$} at 0 -6
\put {$z^{-1}M_N$} at 14 -6
\put {$zN$} at -21 -6
\put {$N$} at -7 -6
\put {$z^{-1}N$} at 7 -6 
\put {$z^{-2}N$} at 21 -6 
\putrule from -21 -10 to 21 -10
\put {$-1$} at -21 -11
\put {0} at -7 -11
\put {1} at 7 -11
\put {2} at 21 -11
\put {$\bullet$} at -21 -10
\put {$\bullet$} at -7 -10
\put {$\bullet$} at 7 -10
\put {$\bullet$} at 21 -10
\put {$\xymatrix@R+3pt {\ar[d]^-{\displaystyle f_N} & \\ &}$} at 1.7 -8
\put {$\xymatrix@R+3pt {\ar[d]^-{\displaystyle z^{-1}f_N} & \\ &}$} at 15.7 -8
\put {$\xymatrix@R+3pt {\ar[d]^-{\displaystyle zf_N} & \\ &}$} at -12.3 -8
\put {$\xymatrix@C+6pt { & \ar[l]_-{\displaystyle z}}$} at 0 -2
\endpicture$$
\vskip2pt

\noindent If $0 \in S^1$ is a regular value of $f$ the Morse function 
$$f_N~=~\overline{f}\vert~:~(M_N;N,z^{-1}N)~=~
\overline{f}^{\,-1}(I;\{0\},\{1\}) \to (I;\{0\},\{1\})$$
has $c_i(f_N)=c_i(f)$ critical points of index $i$.
Each critical point $p \in M$ of $f$ can be regarded as a critical point 
$\overline{p} \in M_N$ of $f_N$. 
For $v \in \GT(f)$ the Morse-Smale complex of $f_N$ with respect 
to $v_N\in \GT(f_N)$ is the based f.\,g. free $\Z$-module chain complex 
$$C^{MS}(M_N,f_N,v_N)~=~F$$ 
with
$$d_F~:~F_i~=~\Z^{c_i(f)} \to F_{i-1}~=~\Z^{c_{i-1}(f)}~;~
\overline{p} \mapsto \sum\limits_{\overline{q}} n(\overline{p},\overline{q})$$
where $n(\overline{p},\overline{q})$ is the algebraic number of 
$v_N$-gradient flow lines in $M_N$ from a critical point 
$\overline{p} \in M_N$ of index $i$ to a 
critical point $\overline{q} \in M_N$ of index $i-1$.

The Novikov complex $C^{Nov}(M,f,v)=C$ counts the
$\overline{v}$-gradient flow lines which start at a critical point
$\overline{p} \in z^jM_N \subset \overline{M}$ and terminate at a
critical point $\overline{q} \in z^kM_N \subset \overline{M}$ with $k
\leq j$.  The Morse-Smale complex $C^{MS}(M_N,f_N,v_N)=F$ only counts
such flow lines with $j=k=0$.  In order to construct $C$ from $F$ we
glue together an infinite number of copies of an 'algebraic fundamental
domain' which gives an algebraic picture of the way the $v_N$-flow
lines enter $M_N$ at $z^{-1}N$ and either die at a critical point of
$f_N$ in $M_N$ or exit at $N$.

As above, given an arbitrary $CW$ structure on $N=f^{-1}(0)$ with
$c_i(N)$ $i$-cells use the handlebody decomposition
$$M_N~=~N \times I \cup \bigcup\limits^m_{i=0}\bigcup_{c_i(f)}D^i\times D^{m-i}~,$$
to give $M_N$ a $CW$ structure with $c_i(N)+c_i(f)$
$i$-cells, one for each $i$-cell of $N$ and one for each
$i$-handle of $(M_N;N,z^{-1}N)$. 
The cellular chain complex of $M_N$ is of the form $C(M_N)=E$ with
$$d_E~=~\pmatrix d_D & c \cr 0 & d_F \endpmatrix~:~
E_i~=~D_i \oplus F_i \to E_{i-1}~=~D_{i-1} \oplus F_{i-1}~,$$
where 
$$\aligned
&D~=~C(N)~~,~~F~=~C(M_N,N\times I)~=~C^{MS}(M_N,f_N,v_N)~,\cr
&D_i~=~\Z^{c_i(N)}~~,~~E_i~=~\Z^{c_i(N)+c_i(f)}~~,~~F_i~=~\Z^{c_i(f)}~,\cr
&g~=~\pmatrix 1 \cr 0\endpmatrix~:~D_i\to E_i~=~D_i \oplus F_i~.
\endaligned$$
Let $h:\alpha D \to E$ be a chain approximation for the
inclusion $h:z^{-1}N \to M$, with components
$$h~=~\pmatrix h_D \cr h_F\endpmatrix~:~\alpha D_i\to E_i~=~D_i \oplus F_i~.$$
The cellular $\Z[z,z^{-1}]$-module chain complex of $\overline{M}$ is 
the algebraic mapping cone
$$C(\overline{M})~=~{\Cal C}(g-zh:D[z,z^{-1}] \to E[z,z^{-1}])~.$$
The chain complex of Farber and Ranicki \cite{FR} was defined to be 
the based f.\,g. free $\Sigma^{-1}\Z[z,z^{-1}]$-module chain complex 
$$C^{FR}(M,f,v,h)~=~\widehat{F}$$
with
$$d_{\widehat{F}}~=~d_F+zh_F(1-zh_D)^{-1}c~:~
\widehat{F}_i~=~\Sigma^{-1}F_i[z,z^{-1}] \to 
\widehat{F}_{i-1}~=~\Sigma^{-1}F_{i-1}[z,z^{-1}]$$
a deformation of $\Sigma^{-1}F[z,z^{-1}]$. 
The inverse of $1-zh_D$ is defined over the Novikov ring $\Z((z))$ by
$$(1-zh_D)^{-1}~=~1+zh_D+z^2(h_D)^2+z^3(h_D)^3+\dots$$
so the induced $\Z((z))$-module chain complex is
$$C^{FR}\big(M,f,v,h;\Z((z))\big)~=~
\Z((z))\otimes_{\Sigma^{-1}\Z[z,z^{-1}]}C^{FR}(M,f,v,h)~=~\overline{F}$$
with
$$d_{\overline{F}}~=~d_F+\sum\limits^{\infty}_{j=1}z^jh_F(h_D)^{j-1}c~:~
\overline{F}_i~=~F_i((z))\to \overline{F}_{i-1}~=~F_{i-1}((z))~.$$
The chain approximation $h$ is an algebraic model for the 
$\overline{v}$-gradient flow across a fundamental domain $(M_N;N,z^{-1}N)$
of $\overline{M}$. The formula for $d_{\widehat{F}}$ is 
interpreted in \S5 as the generating function for the number of 
flow lines in $\overline{M}$ of prescribed length,
with $h_F(h_D)^{j-1}c$ counting the flow lines which start in $M_N$
and terminate in $z^jM_N$, crossing the $j$ walls $N,zN,\dots,z^{j-1}N$
between the adjacent fundamental domains $M_N,zM_N,\dots,z^jM_N$.

The algebraic treatment in \S1 of cokernels of chain maps will
be used in \S6 to prove that the inclusions 
$\widehat{F}_i \to \Sigma^{-1}E_i[z,z^{-1}]$ 
induce the isomorphism of the Cokernel Theorem 6.6
$$C^{FR}(M,f,v,h)~\cong~\hbox{\rm coker}(\phi)$$
with
$$\phi~=~g-zh:\Sigma^{-1}D[z,z^{-1}] \to \Sigma^{-1}E[z,z^{-1}]~.$$
A chain homotopy $k:h \simeq h':\alpha D \to E$ determines an
isomorphism of $\Sigma^{-1}\Z[z,z^{-1}]$-module chain complexes
$$\hbox{\rm coker}(\phi)~\cong~\hbox{\rm coker}(\phi')$$
(Proposition 5.3) giving the isomorphism of the Invariance Theorem 6.7
$$C^{FR}(M,f,v,h)~\cong~\hbox{\rm coker}(\phi)~\cong~
\hbox{\rm coker}(\phi')~\cong~C^{FR}(M,f,v,h')$$
in the special case $N=N'$. The general case is proved by an algebraic
treatment of handle exchanges (Proposition 5.4).

For a gradient-like chain approximation
$$h^{gra}~=~\pmatrix h^{gra}_D \cr h^{gra}_F\endpmatrix~:~\alpha D_i\to
E_i~=~D_i \oplus F_i$$
the algebraic numbers of $\overline{v}$-gradient flow lines between critical
points of $\overline{f}:\overline{M} \to \R$ are given by
$$n(\overline{p},z^j\overline{q})~=~
\cases
\hbox{\rm $(\overline{p},\overline{q})$-coefficient of $d_F:F_i \to F_{i-1}$}
&\hbox{\rm if $j=0$} \cr
\hbox{\rm $(\overline{p},\overline{q})$-coefficient of 
$h^{gra}_F(h^{gra}_D)^{j-1}c:F_i \to F_{i-1}$}
&\hbox{\rm if $j>0$} \cr
0&\hbox{\rm if $j<0$}
\endcases$$
for any critical points $\overline{p},\overline{q} \in M_N$ of $f_N$ 
with index $i,i-1$ respectively.
It follows that for any $-\infty < j <k < \infty$ the Morse-Smale complex of
the real-valued Morse function
$$\aligned
&f_N[j,k]~=~\overline{f}\vert~:\cr
&(M_N[j,k];z^{-j}N,z^{-k}N)~=~\overline{f}^{\,-1}([j,k];\{j\},\{k\})~
=~(\bigcup^{-j}_{\ell=-k}z^{\ell} M_N;z^{-j}N,z^{-k}N)\cr
&\hskip150pt \to ([j,k];\{j\},\{k\})
\endaligned$$
is 
$$\aligned
&C^{MS}(M_N[j,k],f_N[j,k],v_N[j,k])\cr
&\hskip25pt =~(\sum\limits^{-j}_{\ell=-k}z^{\ell}F_i,
d_F+zh^{gra}_F(1+zh^{gra}_D+\dots+z^{k-j-1}(h^{gra}_D)^{k-j-1})c)\cr
&\hskip25pt \cong~\hbox{\rm coker}\big(g-zh^{gra}:
\sum\limits^{-j}_{\ell={-k+1}}z^{\ell}D\to \sum\limits^{-j}_{\ell=-k}z^{\ell}E\big)~.
\endaligned$$
Passing to the direct limit as $k \to \infty$
gives the Morse-Smale complex of the proper real-valued Morse function
$$\aligned
f_N[j,\infty)~=~\overline{f}\vert~:~(M_N[j,\infty),\partial M_N[j,\infty))~
&=~\overline{f}^{\,-1}([j,\infty),\{j\})\cr
&=~(\bigcup^{-j}_{\ell=-\infty}z^{\ell} M_N,z^{-j}N) \to ([j,\infty),\{j\})
\endaligned$$
to be
$$\aligned
C^{MS}(M_N[j,\infty),f_N[j,\infty),v_N[j,\infty))~
&=~\varinjlim_k C^{MS}(M_N[j,k],f_N[j,k],v_N[j,k])\cr
&\cong~\hbox{\rm coker}\big(g-zh^{gra}:
\sum\limits^{-j}_{\ell=-\infty}z^{\ell}D\to 
\sum\limits^{-j}_{\ell=\infty}z^{\ell}E\big)~.
\endaligned$$
Passing to the inverse limit as $j \to -\infty$ gives the Isomorphism
Theorem 6.10 for the Novikov complex, with a basis-preserving
$\Z((z))$-module isomorphism
$$\aligned
C^{Nov}(M,f,v)~
&=~\varprojlim_jC^{MS}\big(M_N[j,\infty),f_N[j,\infty),v_N[j,\infty)\big)\cr
&\cong~\hbox{\rm coker}\big(g-zh^{gra}:D((z)) \to E((z))\big)\cr
&\cong~C^{FR}\big(M,f,v,h^{gra};\Z((z))\big)~.
\endaligned
$$
The geometric differential in $C=C^{Nov}(M,f,v)$ is just the algebraic 
differential in 
$\overline{F}=C^{FR}\big(M,f,v,h^{gra};\Z((z)))$, with
$$\aligned
d_C(\overline{p})~
&=~\sum\limits_{\overline{q}}\sum\limits_{j\in \Z}
 n(\overline{p},z^j\overline{q})z^j\overline{q}\cr
&=~(d_F+\sum\limits^{\infty}_{j=1}z^jh^{gra}_F(h^{gra}_D)^{j-1}c)(\overline{p})\cr
&=~(d_F+zh^{gra}_F(1-zh^{gra}_D)^{-1}c)(\overline{p})\cr
&=~d_{\overline{F}}(\overline{p}) \in \overline{F}_{i-1}~=~\Z^{c_{i-1}(f)}~,
\endaligned$$
so that there is also a basis-preserving $\Sigma^{-1}\Z[z,z^{-1}]$-module isomorphism
$$C^{Paj}(M,f,v)~\cong~C^{FR}(M,f,v,h^{gra})~.$$
The projection
$$\aligned
p~:~C\big(\overline{M};\Z((z))\big)~
&=~{\Cal C}\big(g-zh^{gra}:D((z)) \to E((z))\big)\cr
&\to \hbox{\rm coker}\big(g-zh^{gra}:D((z)) \to E((z))\big)~=~C^{Nov}(M,f,v)
\endaligned$$
pieces together $\overline{v}$-gradient flow lines in
$\overline{M}$ from their intersections with the translates
$z^jM_N\subset \overline{M}$ ($j \in \Z$) of the fundamental domain
$M_N$\,; $p$ is a chain equivalence with torsion
$$\aligned
\tau(p)~
&=~\sum\limits^{\infty}_{i=0}(-)^{i+1}\tau\big(1-zh^{gra}_D:D_i((z))
\to D_i((z))\big)\cr
&=~\prod\limits^{\infty}_{i=0}\text{det}\big(1-zh^{gra}_D:D_i((z)) \to
D_i((z))\big)^{(-)^{i+1}}\cr
& \in \widehat{W}(\Z) \subset K_1\big(\Z((z))\big)~=~K_1(\Z)\oplus
K_0(\Z) \oplus \widehat{W}(\Z)
\endaligned$$
with $\widehat{W}(\Z)=1+z\Z[[z]]$ under multiplication.
The kernel of $p$ corresponds to the closed orbits of the
$v$-gradient flow lines in $M$, which avoid the critical points of
$f$ and so do not contribute to the Novikov complex, and which are
counted by the torsion of $p$.

I am grateful to Andrei Pajitnov for valuable suggestions for improving
the preprint version of the paper.

\subheading{1. Cones and cokernels}

The Novikov complex of a circle-valued Morse function $f:M \to S^1$ will be
shown in \S6 to be isomorphic to the cokernel of a chain map
constructed from the Morse-Smale complex of a fundamental domain
for the infinite cyclic cover $\overline{M}=f^*\R$ of $M$;
the algebraic mapping cone of the chain map is a cellular chain complex
of $\overline{M}$.
This section is accordingly devoted to the relationship
between the algebraic mapping cone and the cokernel of a chain map.
The algebraic mapping cone is a chain homotopy invariant.
The cokernel is not a chain homotopy invariant, although it is a
homology invariant.
The main novelty of this section is the introduction of an 
equivalence relation on chain maps called `chain isotopy', 
which is stronger than chain homotopy, and is
such that the cokernels of chain isotopic maps are isomorphic.
The chain map with cokernel the Novikov complex depends on a choice of 
chain map in a chain homotopy class\,; a different choice will give a chain isotopic 
chain map, with isomorphic cokernel.

Let $A$ be a ring. The {\it algebraic mapping cone} of an $A$-module chain map 
$$\phi :D \to E$$
is the $A$-module chain complex ${\Cal C}(\phi)$ defined by
$$d_{{\Cal C}(\phi)}~=~\pmatrix d_E & (-)^{i-1}\phi  \cr 0 & d_D \endpmatrix
~:~{\Cal C}(\phi)_i~=~E_i \oplus D_{i-1} \to {\Cal C}(\phi)_{i-1}~=~
E_{i-1} \oplus D_{i-2}~.$$
The natural projections
$$p~:~{\Cal C}(\phi)_i~=~E_i \oplus D_{i-1} \to 
\hbox{\rm coker}(\phi :D_i \to E_i)~;~(x,y) \mapsto [x]$$
define a chain map
$$p~:~{\Cal C}(\phi) \to \hbox{\rm coker}(\phi)~.$$

\proclaim{\bf Proposition 1.1} 
For an injective chain map $\phi:D \to E$ the natural projection 
$p:{\Cal C}(\phi) \to \hbox{\rm coker}(\phi)$
induces isomorphisms in homology
$$p_*~:~H_*\big({\Cal C}(\phi)\big)~ \cong~ H_*\big(\hbox{\rm coker}(\phi)\big)~.$$
\endproclaim
\demo{Proof} The short exact sequences of $A$-module chain complexes
$$\aligned
&0 \to E \to {\Cal C}(\phi) \to D_{*-1} \to 0\cr
&0 \to D @>\phi >> E \to \hbox{\rm coker}(\phi) \to 0 
\endaligned$$
induce long exact sequences in homology, which are related by a
natural transformation
$$\xymatrix@C-10pt
{\dots \ar[r] & H_i(D) \ar[d]^1 \ar[r]^{\phi _*} &
H_i(E) \ar[d]^1 \ar[r] &
H_i({\Cal C}(\phi)) \ar[d]^{p_*} \ar[r] & H_{i-1}(D) \ar[d]^1 \ar[r] & \dots \\
\dots \ar[r] & H_i(D)  \ar[r]^{\phi _*} &
H_i(E)  \ar[r] &
H_i\big(\hbox{\rm coker}(\phi)\big) \ar[r] & H_{i-1}(D) \ar[r] & \dots}$$
It now follows from the 5-lemma that the induced morphisms
$p_*$ are isomorphisms.\qed\enddemo

As usual, a {\it chain homotopy} between
chain maps $\phi,\phi':D \to E$ 
$$\theta~:~\phi~\simeq~\phi'~:~D \to E$$
is a collection of $A$-module morphisms $\theta:D_i \to E_{i+1}$ such that
for each $i$
$$\phi' - \phi~=~d_E\theta+\theta d_D ~:~ D_i \to E_i~.$$
Chain homotopic chain maps have isomorphic algebraic mapping cones\, :

\proclaim{\bf Proposition 1.2} A chain homotopy
$\theta:\phi \simeq \phi':D \to E$
determines an isomorphism of the algebraic mapping cones
$$I~=~\pmatrix 1 & \pm \theta  \cr 0 & 1 \endpmatrix~:~
{\Cal C}(\phi) \to {\Cal C}(\phi')~.$$
If $D,E$ are based f.\,g. free then $I$ is a simple isomorphism.
\endproclaim
\demo{Proof} By construction. \qed\enddemo

In general, the cokernels of chain homotopic chain maps are not isomorphic
(or even chain equivalent). The following relation will be convenient in 
dealing with cokernels of chain maps, in order to avoid this problem.

\noindent{\bf Definition 1.3}
A {\it chain isotopy} between chain maps $\phi,\phi':D \to E$
$$\psi~:~\phi \sim \phi'~:~D \to E$$ 
is a collection of $A$-module morphisms $\psi:E_i \to E_{i+1}$ such that
\roster
\item"(i)" for each $i$
$$\phi'~=~(1+d_E\psi+\psi d_E)\phi~ :~ D_i \to E_i~,$$ 
defining a chain homotopy
$$\psi\phi~ :~\phi~ \simeq~ \phi'~:~D \to E~,$$
\item"(ii)" each
$$1+d_E\psi+\psi d_E~:~E_i \to E_i$$ 
is an automorphism.\hfill\qed
\endroster

\proclaim{\bf Proposition 1.4}
Chain isotopy is an equivalence relation on chain maps.
\endproclaim
\demo{Proof} Reflexivity: every chain map $\phi:D \to E$ is
isotopic to itself by $0:\phi \sim \phi$.\newline
Symmetry: for any chain isotopy $\psi:\phi\sim\phi':D \to E$ 
define a chain isotopy $\psi^-:\phi'\sim \phi$ by
$$\psi^-~=~-(1+d_E\psi+\psi d_E)^{-1}\psi~:~E_i \to E_{i+1}~,$$
with 
$$1+d_E\psi^-+\psi^- d_E~=~(1+d_E\psi+\psi d_E)^{-1}~:~E_i \to E_i~.$$ 
Transitivity: for any chain isotopies 
$\psi:\phi\sim\phi'$, $\psi':\phi'\sim\phi'':D \to E$ 
define a chain isotopy $\psi'':\phi\sim \phi''$ by
$$\psi''~=~\psi+\psi'(1+d_E\psi+\psi d_E)~:~E_i \to E_{i+1}$$
with
$$1+d_E\psi''+\psi'' d_E~=~(1+d_E\psi'+\psi' d_E)(1+d_E\psi+\psi d_E)
~:~E_i \to E_i~.\qed$$ 
\enddemo

Isotopic chain maps have isomorphic cokernels\, :

\proclaim{\bf Proposition 1.5} 
A chain isotopy $\psi:\phi\sim\phi':D \to E$ determines isomorphisms 
of chain complexes
$$\aligned
&q~=~\pmatrix 1 & \pm \psi\phi  \cr 0 & 1 \endpmatrix~:~
{\Cal C}(\phi) \to {\Cal C}(\phi')~,\cr
&r~=~[1+d_E\psi+\psi d_E]~:~
\hbox{\rm coker}(\phi) \to \hbox{\rm coker}(\phi')
\endaligned$$
and a chain homotopy
$$s~:~rp~\simeq~p'q~:~{\Cal C}(\phi) \to \hbox{\rm coker}(\phi')$$
with $p:{\Cal C}(\phi) \to \hbox{\rm coker}(\phi)$,
$p':{\Cal C}(\phi') \to \hbox{\rm coker}(\phi')$ the projections.
\endproclaim
\demo{Proof} 
The isomorphism $q$ is a special case of 1.2.\newline
The isomorphism $r$ is given by the morphism of exact sequences
$$\xymatrix@C+40pt @R+5pt{
D \ar@{=}[d] \ar[r]^{\displaystyle{\phi}} & 
E \ar[d]_{\displaystyle{\cong}}^-{\displaystyle{1+d_E\psi +\psi d_E}} \ar[r] 
&\hbox{\rm coker}(\phi)
\ar[d]^-{\displaystyle{r}}_-{\displaystyle{\cong}} \\
D \ar[r]^{\displaystyle{\phi'}} & E \ar[r] 
&\hbox{\rm coker}(\phi')}\hskip-5pt
\xymatrix@R+7pt {\ar[r] & 0 \\ \ar[r] & 0}$$
The $A$-module morphisms
$$s~:~{\Cal C}(\phi)_i~=~E_i \oplus D_{i-1} \to \hbox{\rm coker}(\phi':
D_{i+1} \to E_{i+1})~;~(x,y) \mapsto [\psi(x)]$$
define a chain homotopy
$$s~:~rp ~\simeq~p'q~:~{\Cal C}(\phi) \to \hbox{\rm coker}(\phi')$$
in the diagram
$$\xymatrix@C+10pt @R+10pt{
{\Cal C}(\phi) \ar[d]_-{\displaystyle{q}} \ar[r]^-{\displaystyle{p}}
\ar@{~>}[dr]^{\displaystyle{s}}
& \hbox{\rm coker}(\phi) \ar[d]^-{\displaystyle{r}}\\
{\Cal C}(\phi') \ar[r]^-{\displaystyle{p'}} & \hbox{\rm coker}(\phi')~.}
\eqno{\lower52pt\hbox{\qed}}$$
\enddemo

\noindent{\bf Definition 1.6} An {\it embedding} of chain complexes
is a chain map $\phi:D \to E$ such that each $\phi:D_i \to E_i$ 
is a split injection.\qed

\proclaim{\bf Proposition 1.7} Let $\phi:D \to E$ be an embedding of 
$A$-module chain complexes.\newline
{\rm (i)} The natural projection $p:{\Cal C}(\phi) \to \hbox{\rm coker}(\phi)$ is a 
chain equivalence.\newline
{\rm (ii)} If each $D_i,E_i,\hbox{\rm coker}(\phi :D_i \to E_i)$ is a based f.\,g. free
$A$-module, then $p$ is a chain equivalence with torsion
$$\tau(p)~=~\sum\limits^{\infty}_{i=0}(-)^i\tau({\Cal E}_i) \in K_1(A)$$
where $\tau({\Cal E}_i)\in K_1(A)$ is the torsion of the short exact sequence of
based f.\,g. free $A$-modules
$$\xymatrix{
{\Cal E}_i~:~0 \ar[r] & D_i \ar[r]^-{\displaystyle{\phi}} &
E_i \ar[r] & \hbox{\rm coker}(\phi:D_i \to E_i) \ar[r] & 0~.}$$
{\rm (iii)} If $\psi:\phi \sim \phi':D \to E$ is a chain isotopy 
then $\phi':D \to E$ is also an embedding.
With bases as in {\rm (ii)}, and the isomorphism given by 1.5
$$r~=~[1+d_E\psi+\psi d_E]~:~\hbox{\rm coker}(\phi)~
\cong~\hbox{\rm coker}(\phi')$$
has torsion
$$\tau(r)~=~\sum\limits^{\infty}_{i=0}(-)^i
(\tau({\Cal E}'_i)-\tau({\Cal E}_i) )\in K_1(A)~.$$
\endproclaim
\demo{Proof}  (i) Extend each $\phi :D_i \to E_i$ to a direct sum system
$$D_i \vcenter{\hbox{$\phi \atop \ra{2}$}
\hbox{$\la{2} \atop e$}}~E_i~
\vcenter{\hbox{$j\atop \ra{2}$}
\hbox{$\la{2} \atop k$}}~F_i$$
with $F_i=\hbox{\rm coker}(\phi :D_i \to E_i)$. Let
$$c~=~ed_Ek~:~F_i \to D_{i-1}~~,~~d_F~=~jd_Ek~:~F_i \to F_{i-1}~,$$
so that there is defined an isomorphism of chain complexes 
$(\phi ~k):E' \to E$ with
$$d_{E'}~=~\pmatrix d_D & c  \cr 0 & d_F \endpmatrix~:~
E'_i~=~D_i \oplus F_i \to E'_{i-1}~=~D_{i-1} \oplus F_{i-1}~.$$
The chain map $p:{\Cal C}(\phi) \to F=\hbox{\rm coker}(\phi)$ is given by
$$p~=~\left( \matrix j & 0 \endmatrix \right) ~:~{\Cal C}(\phi)_i~=~
E_i \oplus D_{i-1} \to F_i~.$$
The chain map $g:F \to {\Cal C}(\phi)$ defined by
$$g~=~\pmatrix k  \cr -c \endpmatrix~:~F_i \to 
{\Cal C}(\phi)_i ~=~E_i \oplus D_{i-1}$$
is such that 
$$pg~=~1~:~F \to F~~,~~h~:~g p ~\simeq~ 1~:~{\Cal C}(\phi) \to {\Cal C}(\phi)$$
with $h$ the chain homotopy
$$h~=~\pmatrix 0  & 0  \cr e & 0 \endpmatrix~:~
{\Cal C}(\phi)_i ~=~E_i \oplus D_{i-1}\to {\Cal C}(\phi)_{i+1} ~=~
E_{i+1} \oplus D_i.$$
Thus $p:{\Cal C}(\phi) \to F$, $g:F \to {\Cal C}(\phi)$ are inverse chain equivalences.\newline
(ii) Immediate from the 0-dimensional case, which follows from the sum
formula $\tau(fg)=\tau(f)+\tau(g)$.\newline
(iii) By the chain homotopy invariance of torsion, the sum formula,
the identity $\tau(q)=0$ with $q$ as in 1.5 and (ii)
$$\aligned
\tau(r)~&=~\tau(p') - \tau(p) + \tau(q)\cr
&=~ \tau(p')-\tau(p)\cr
&=~ \sum\limits^{\infty}_{i=0}(-)^i\tau({\Cal E}'_i) 
- \sum\limits^{\infty}_{i=0}(-)^i\tau({\Cal E}_i) \in K_1(A)~.\hfill\qed
\endaligned$$
\enddemo

\proclaim{\bf Proposition 1.8} Let $E$ be an $A$-module chain complex of
the form
$$d_E~=~\pmatrix d_D & c \cr 0 & d_D \endpmatrix~:~
E_i~=~D_i \oplus F_i \to E_{i-1}~=~D_{i-1} \oplus F_{i-1}~.$$
Let $\phi:D \to E$ be a chain map, with
$$\phi~=~\pmatrix \phi_D \cr \phi_F \endpmatrix~:~
D_i \to E_i~=~D_i \oplus F_i~.$$
{\rm (i)} If each $\phi_D:D_i \to D_i$ is an automorphism then 
\roster
\item"(a)" $\phi$ is an embedding of chain complexes.
\item"(b)" The chain complex $\widehat{F}$ defined by
$$d_{\widehat{F}}~=~d_F - \phi_F(\phi_D)^{-1}c ~:~\widehat{F}_i~=~F_i \to \widehat{F}_{i-1}~=~F_{i-1}$$
is such that the inclusions $F_i \to E_i$ induce an isomorphism of chain complexes
$$\widehat{F}~ \cong~ \hbox{\rm coker}(\phi)~.$$
\item"(c)" The natural projection $p:{\Cal C}(\phi) \to \hbox{\rm coker}(\phi)$ 
is a chain equivalence. The chain complex $K$ defined by
$$\aligned
&d_K~=~\pmatrix d_D & 0 \cr (-1)^i\phi_D & \phi_Dd_D(\phi_D)^{-1} \endpmatrix~:\cr
&\hskip50pt K_i~=~D_{i-1} \oplus D_i \to K_{i-1}~=~D_{i-2} \oplus D_{i-1}
\endaligned$$
{\rm (}i.e. the algebraic mapping cone of the isomorphism of chain complexes
$\phi_D:D=(D_i,d_D) \to (D_i,\phi_Dd_D(\phi_D)^{-1})${\rm )}
is contractible, and fits into an exact sequence
$$0 \to K \to {\Cal C}(\phi) \xymatrix{\ar[r]^{\displaystyle{p}}&} 
\hbox{\rm coker}(\phi) \to 0$$
with
$$\aligned
&K_i~=~D_{i-1} \oplus D_i \to {\Cal C}(\phi)_i~=~D_{i-1} \oplus D_i \oplus E_i~;\cr
&\hskip130pt (x,y) \mapsto (x,y,\phi_F(\phi_D)^{-1}(y))~.
\endaligned$$
\item"(d)" If $D,F$ are based f.\,g. free the natural projection
$$p~:~{\Cal C}(\phi) \to \hbox{\rm coker}(\phi)~\cong~\widehat{F}$$ 
is a chain equivalence of based f.\,g. free $A$-module chain complexes 
with torsion
$$\tau(p:{\Cal C}(\phi) \to \widehat{F})~=~
-\tau(K)~=~-\sum\limits^{\infty}_{i=0}(-)^i\tau(\phi_D:D_i \to D_i) \in K_1(A)~.$$
\endroster
{\rm (ii)} Given a chain homotopy of chain maps
$$\theta~:~\phi~\simeq~ \phi'~:~D \to E$$ 
with $\phi_D,\phi'_D:D_i \to D_i$ automorphisms write
$$\theta~=~\pmatrix \theta_D \cr \theta_F \endpmatrix~:~
D_i \to E_{i+1}~=~D_{i+1} \oplus F_{i+1}~.$$
The morphisms defined by
$$\psi~=~\pmatrix \theta_D(\phi_D)^{-1} & 0 \cr 
\theta_F(\phi_D)^{-1} & 0 \endpmatrix~:~E_i~=~D_i \oplus F_i \to
E_i~=~D_i \oplus F_i$$
are such that
$$\phi'~=~(1+d_E\psi + \psi d_E)\phi~:~D_i \to E_i~.$$
Thus if each $1+d_E\psi+\psi d_E:E_i \to E_i$ is an automorphism there is
defined a chain isotopy of embeddings
$$\psi~:~\phi~\sim~\phi'~:~D \to E$$
and as in 1.5 there is defined an isomorphism
$$r~=~[1+d_E\psi+\psi d_E]~:~\hbox{\rm coker}(\phi)~\cong~\hbox{\rm coker}(\phi')~.$$
Moreover, if $D,F$ are based f.g. free the isomorphism
$$r~:~\widehat{F}~\cong~\hbox{\rm coker}(\phi)~\cong~\hbox{\rm coker}(\phi')~
\cong~\widehat{F}'$$
is simple, with
$$\aligned
\tau(r)~&=~\sum\limits^{\infty}_{i=0}(-)^i
(\tau(\phi_D:D_i \to D_i)-\tau(\phi'_D:D'_i \to D'_i))\cr
&=~0  \in K_1(A)~.
\endaligned$$
\endproclaim
\demo{Proof} (i) (a) Each $\phi :D_i \to E_i$ is a split injection.\newline
(b) It is clear that the inclusion $F \to E$ induces a chain map
$\widehat{F} \to \hbox{\rm coker}(\phi)$. 
The $A$-module morphisms
$$\pmatrix -\phi_F(\phi_D)^{-1} &  1\endpmatrix~:~E_i~=~D_i \oplus F_i \to \widehat{F}_i~=~F_i$$
induce the inverse chain isomorphism
$\hbox{\rm coker}(\phi) \to \widehat{F}$.\newline
(c) Immediate from 1.7 (i).\newline
(d) Apply 1.7 (ii), noting that the short exact sequence
$$\xymatrix{{\Cal E}_i~:~0 \ar[r] & }
\xymatrix@C+20pt{
D_i \ar[r]^-{\displaystyle{\pmatrix \phi_D \cr \phi_F \endpmatrix}} &}
\xymatrix@C+45pt{D_i \oplus F_i
\ar[r]^-{\displaystyle{(-\phi_F(\phi_D)^{-1}~1)}} &}
\xymatrix{F_i \ar[r] & 0}$$
has torsion
$$\tau({\Cal E}_i)~=~-\tau(\phi_D:D_i \to D_i) \in K_1(A)~.$$
(ii) It follows from
$$\aligned
&\pmatrix d_D & c \cr 0 & d_F \endpmatrix 
\pmatrix \phi_D \cr \phi_E \endpmatrix~=~\pmatrix \phi_D \cr \phi_E \endpmatrix d_D~:~
D_i \to E_{i-1}~=~D_{i-1} \oplus F_{i-1}~,\cr
&\pmatrix \phi'_D \cr \phi'_F \endpmatrix -
\pmatrix \phi_D \cr \phi_F \endpmatrix~=~
\pmatrix d_D & c \cr 0 & d_F \endpmatrix 
\pmatrix \theta_D \cr \theta_E \endpmatrix +\pmatrix \theta_D \cr 
\theta_E \endpmatrix d_D~:~D_i \to E_i~=~D_i \oplus F_i
\endaligned$$
that
$$\aligned
&(1+d_E\psi+\psi d_E)\phi~=~
\bigg( \pmatrix 1 & 0 \cr 0 & 1 \endpmatrix +
\pmatrix d_D & c \cr 0 & d_F \endpmatrix 
\pmatrix \theta_D(\phi_D)^{-1} & 0 \cr 
\theta_F(\phi_D)^{-1} & 0 \endpmatrix\cr
&\hskip150pt  +\pmatrix \theta_D(\phi_D)^{-1} & 0 \cr 
\theta_F(\phi_D)^{-1} & 0 \endpmatrix
\pmatrix d_D & c \cr 0 & d_F \endpmatrix \bigg)\pmatrix \phi_D \cr \phi_E \endpmatrix
\cr
&\hphantom{(1+d_E\psi+\psi d_E)\phi~}
=~ \pmatrix \phi'_D \cr \phi'_F \endpmatrix ~=~
\phi'~:~D_i \to E_i~=~D_i \oplus F_i~.
\endaligned$$
If $1+d_E\psi + \psi d_E:E_i \to E_i$ is an automorphism then
$\psi:\phi\simeq \phi':D \to E$ is a chain isotopy, and by 1.7 (iii)
$$\aligned
\tau(r:\widehat{F} \cong \widehat{F}')~&=~
\sum\limits^{\infty}_{i=0}(-)^i(\tau({\Cal E}'_i)-\tau({\Cal E}_i) )\cr
&=~\sum\limits^{\infty}_{i=0}(-)^i(\tau(\phi_D:D_i \to D_i)
-\tau(\phi'_D:D'_i \to D'_i)) \in K_1(A)~.
\endaligned$$
Now there is defined an isomorphism of short exact sequences
$$\xymatrix@R+15pt{{\Cal E}_i~:~0 \ar[r] & \\
{\Cal E}'_i~:~0 \ar[r]& }
\xymatrix@C+20pt @R+15pt{
D_i \ar@{=}[d]
\ar[r]^-{\displaystyle{\pmatrix \phi_D \cr \phi_F \endpmatrix}} &\\
D_i \ar[r]^-{\displaystyle{\pmatrix \phi'_D \cr \phi'_F \endpmatrix}} &}
\xymatrix@C+45pt @R+15pt{
D_i \oplus F_i \ar[d]^-{\displaystyle{1+d_E\psi +\psi d_E}}
\ar[r]^-{\displaystyle{(-\phi_F(\phi_D)^{-1}~~1)}} &\\
D_i \oplus F_i 
\ar[r]^-{\displaystyle{(-\phi'_F(\phi'_D)^{-1}~~1)}} &}
\xymatrix@R+12pt{\hbox{$\widehat{F}$}_i \ar[d]^-{\displaystyle{r}} \ar[r] & 0\\
\hbox{$\widehat{F}'$}_i \ar[r] & 0}$$
so that
$$\tau(r:\widehat{F}_i \to \widehat{F}'_i)~=~
\tau(1+d_E\psi+\psi d_E:E_i \to E_i) \in K_1(A)~.$$
The chain complex automorphism $1+d_E\psi+\psi d_E:E \to E$ is chain
homotopic to $1:E \to E$, so
$$\aligned
\tau(r:\widehat{F} \to\widehat{F}')~&=~
\sum\limits^{\infty}_{i=0}(-1)^i\tau(r:\widehat{F}_i \to \widehat{F}'_i)\cr
&=~\sum\limits^{\infty}_{i=0}(-1)^i\tau(1+d_E\psi+\psi d_E:E_i \to E_i)\cr
&=~\tau(1+d_E\psi+\psi d_E:E \to E)\cr
&=~\tau(1:E \to E)~=~0 \in K_1(A)~.
\endaligned\eqno{\lower37pt\hbox{\qed}}$$
\enddemo
The formula $d_{\widehat{F}}=d_F - \phi_F(\phi_D)^{-1}c:F_i \to F_{i-1}$
was first obtained in \cite{FR,2.3}.

\subheading{2. The Morse-Smale complex $C^{MS}(M,f,v)$}

This section recalls the properties of the Morse-Smale complex
$C^{MS}(M,f,v)$ of a real-valued Morse function on an 
$m$-dimensional cobordism
$$f~:~(M;N,N') \to (I;\{0\},\{1\})$$
with respect to any $v \in \GT(f)$.
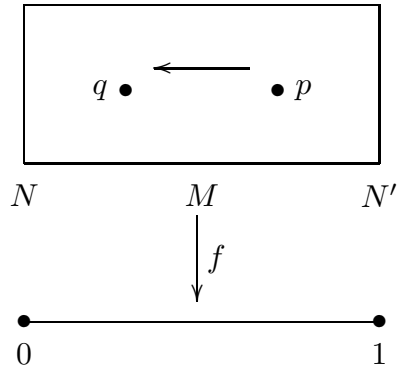
$$\beginpicture
\setcoordinatesystem units <8pt,10pt>  
\putrule from -7 0 to 7 0
\putrule from -7 0 to -7 -5
\putrule from 7 0 to 7 -5
\putrule from -7 -5 to 7 -5
\put {$\bullet$} at -3 -2.7
\put {$q$} at -4 -2.7
\put {$\bullet$} at 3 -2.7
\put {$p$} at 4 -2.7
\put {$M$} at 0 -6
\put {$N$} at -7 -6
\put {$N'$} at 7 -6 
\putrule from -7 -10 to 7 -10
\put {0} at -7 -11
\put {1} at 7 -11
\put {$\bullet$} at -7 -10
\put {$\bullet$} at 7 -10
\put {$\xymatrix@R+3pt {\ar[d]^-{\displaystyle f} & \\ &}$} at 1.7 -8
\put {$\xymatrix@C+6pt {& \ar[l]}$} at 0 -2
\endpicture$$
\vskip2pt

\noindent{\bf Definition 2.1}
The {\it Morse-Smale complex} of $f:M \to \R$ with respect to $v \in \GT(f)$ 
is the based f.\,g. free $\Z[\pi_1(M)]$-module 
chain complex $C^{MS}(M,f,v)$ with
\roster
\item"(i)" $\hbox{\rm rank}_{\Z[\pi_1(M)]}C^{MS}_i(M,f,v)=c_i(f)$, 
with one basis element $\widetilde{p}$ for each critical
point $p \in M$ of index $i$, corresponding to a choice
of lift $\widetilde{p} \in \widetilde{M}$,
\item"(ii)" the boundary $\Z[\pi_1(M)]$-module morphisms 
$$\aligned
d~:~C^{MS}_i(M,f,v)~=~\Z[\pi_1(M)]^{c_i(f)}&\to 
C^{MS}_{i-1}(M,f,v)~=~\Z[\pi_1(M)]^{c_{i-1}(f)}~;\cr
&\widetilde{p} \mapsto \sum\limits_{\widetilde{q}}n(\widetilde{p},\widetilde{q})\widetilde{q}
\endaligned$$
with $n(\widetilde{p},\widetilde{q}) \in \Z$ the algebraic number of 
$\widetilde{v}$-gradient flow lines in $\widetilde{M}$ joining $\widetilde{p}$ to 
$\widetilde{q}$.\qed
\endroster

\medskip

A Morse function on a cobordism
$$f~:~(M;N,N') \to (I;\{0\},\{1\})$$
and $v \in \GT(f)$ determine a handlebody decomposition
$$M~=~N \times I \cup\bigcup\limits^m_{i=0}
\bigcup\limits_{c_i(f)}h^i$$
with $c_i(f)$ $i$-handles $h^i=D^i \times D^{m-i}$.
Given a $CW$ structure on $N$ with $c_i(N)$ $i$-cells $e^i \subset N$ let $M$
have the $CW$ structure with 
$$c_i(M)~=~c_i(N)+c_i(f)$$ 
$i$-cells\,: there is one $i$-cell $e^i \times I \subset M$ for each
$i$-cell $e^i \subset N$, and one $i$-cell $h^i \subset M$ for each
critical point of index $i$.  Let $\widetilde{M}$ be the universal
cover of $M$, and let $\widetilde{N},\widetilde{N}'$ be the
corresponding covers of $N,N'$.

\proclaim{\bf Proposition 2.2} 
The Morse-Smale complex $C^{MS}(M,f,v)$ is 
the relative cellular chain complex of $(\widetilde{M},\widetilde{N}\times I)$
$$C^{MS}(M,f,v)~=~C(\widetilde{M},\widetilde{N}\times I)~.$$
\endproclaim
\demo{Proof} See Franks \cite{Fr}.\hfill\qed
\enddemo

\noindent{\bf Definition 2.3}  Given a Morse function
$f:(M;N,N') \to (I;\{0\},\{1\})$, $v \in \GT(f)$ 
and $CW$ structures on $N$ and $N'$ write
$$D~=~C(\widetilde{N})~~,~~D'~=~C(\widetilde{N}')~,~E~=~C(\widetilde{M})~~,~~
F~=~C^{MS}(M,f,v)~=~C(\widetilde{M},\widetilde{N}\times I)~.$$
(i) The cellular chain complex  of $\widetilde{M}$ is of the form
$$d_E~=~\pmatrix d_D & c \cr 0 & d_F \endpmatrix~:~
E_i~=~D_i \oplus F_i \to E_{i-1}~=~D_{i-1}\oplus F_{i-1}$$
for a {\it birth} chain map $c:F_{*+1} \to D$, that is 
$$E~=~C(\widetilde{M})~=~{\Cal C}(c)~.$$
(ii) The natural map $g:N \to M$ is the inclusion of a $CW$ subcomplex,
inducing the embedding
$$g~=~\pmatrix 1 \cr 0 \endpmatrix~:~
C(\widetilde{N})~=~D \to C(\widetilde{M})~=~E~.$$
A {\it crossing chain approximation}
$$h~:~C(\widetilde{N}')~=~D' \to C(\widetilde{M})~=~E$$
is a chain map induced by the inclusion 
$h:N' \to M$ -- in general this is not the inclusion of a $CW$ 
subcomplex, so the construction requires the cellular approximation theorem,
and the $CW$ structures only determine a chain map $h$ up to chain homotopy. 
The components of 
$$h~=~\pmatrix h_D \cr h_F \endpmatrix~:~D'_i \to E_i~=~D_i \oplus F_i$$
are such that
$$\aligned
&d_Dh_D+ch_F~=~h_Dd_{D'}~:~D'_i \to D_{i-1}~,\cr
&d_Fh_F~=~h_Fd_{D'}~:~D'_i \to F_{i-1}
\endaligned$$
defining a {\it death} chain map
$$h_F~:~D' \to F$$ 
and a {\it survival} chain homotopy
$$h_D~:~ch_F~ \simeq~ 0~:~D' \to D_{*-1}~.\eqno{\qed}$$

\noindent{\bf Remark 2.4} For a Morse function
$$f~:~(M;N,N') \to (I;\{0\},\{1\})$$
with $v \in \GT(f)$ there are 4 types of $v$-gradient flow
lines, which we shall call $\alpha ,\beta,\gamma,\delta$,
corresponding to the 4 morphisms
$$d_F~:~F_i \to F_{i-1}~,~c~:~F_i \to D_{i-1}~,~h_F~:~D'_i \to F_i~,~
h_D~:~D'_i \to D_i~,$$
as follows :
\roster
\item"(a)" complete: $\alpha $ starts at an index $i$ critical point $p\in M$
and terminates at an index $i-1$ critical point $q \in M$. 
\item"(b)" birth: $\beta$ starts at an index $i$ critical point $p \in M$
and terminates in $N$. 
\item"(c)" death: $\gamma$ starts in $N'$ and terminates at an index 
$i$ critical point $p \in M$.
\item"(d)" survival: $\delta$ starts in $N'$ and terminates in $N$.
\endroster
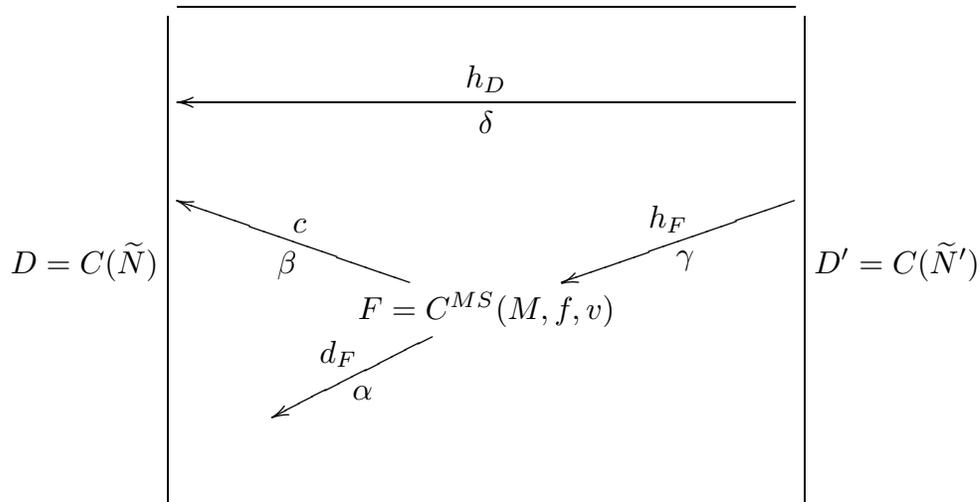
$$\xymatrix
{ \ar@{-}[rrrr]
\ar@{-}[ddddd]_-{\displaystyle{D=C(\widetilde{N})}} & & & &
\ar@{-}[ddddd]^-{\displaystyle{D'=C(\widetilde{N}')}} \\
& & & & \ar[llll]^-{\displaystyle{\delta}}_-{\displaystyle{h_D}}   \\
& & & & \ar[dll]^-{\displaystyle{\gamma}}_-{\displaystyle{h_F}}  \\
& & F=C^{MS}(M,f,v) \ar[llu]^-{\displaystyle{\beta}}_-{\displaystyle{c}}  
\ar[dl]^-{\displaystyle{\alpha }}_-{\displaystyle{d_F}} & & \\
& & & & \\
\ar@{-}[rrrr] & & & & }$$
A crossing chain approximation $h:D' \to E$ corresponds to the flow
lines which start in $N'$, i.e.  those of death and survival type. 
Pajitnov \cite{P4,\S4} obtained an analogue of the cellular
approximation theorem for the gradient flow : for any Morse function
$f:M \to S^1$ with regular value $0 \in S^1$ it is possible to choose
$v \in \GT(f)$ and handlebody structures on $N$ and $N'$ such that
\roster
\item"(i)" every survival $v$-gradient flow line in $M$ which starts in an
$i'$-handle of $N'$ ends in an $i$-handle of $N$ with $i \leq i'$,
\item"(ii)" there is a finite number of rel $\partial$ homology classes
of survival flow lines as in (i) with $i=i'-1$,
\item"(iii)" there exist $d_F,c,h_F$ as in 2.3 which actually count the 
flow lines of type $\alpha ,\beta,\gamma$, and $h_D$ 
counts the rel $\partial$ homology classes of survival flow lines $\delta$
as in (ii),
\item"(iv)" the function which sends $y \in M\backslash \hbox{Crit}(f)$ 
to the endpoint $\Phi(y) \in N \cup \hbox{Crit}(f)$ of the flow line of $f$ 
which starts at $y$
$$\Phi~:~M\backslash \hbox{Crit}(f) \to N \cup \hbox{Crit}(f)~;~y \mapsto \Phi(y)$$
restricts to a function
$$N'\cap \Phi^{-1}(N) \to N \cap \Phi(N')~;~x' \mapsto \Phi(x')$$
which is a partially defined map $N' \to N$ inducing the
'partial chain map' 
$$h_D~:~D'~=~C(\widetilde{N}') \to D~=~C(\widetilde{N})~.$$
\endroster
(See also Hutchings and Lee \cite{HL}).\hfill\qed

We shall now express the Morse-Smale complex $C^{MS}(M,f,v)$ of a Morse function
$$f~:~(M;N,N'') \to ([0,2];0,2)~~(v \in \GT(f))$$
which is transverse regular at $1 \in [0,2]$ in terms of the Morse-Smale 
complexes $C^{MS}(M',f',v')$, $C^{MS}(M'',f'',v'')$ of the restrictions
$$\aligned
&f'~=~f\vert~:~(M';N,N')~=~f^{-1}([0,1];\{0\},\{1\}) \to ([0,1];\{0\},\{1\})~,\cr
&f''~=~f\vert~:~(M'';N',N'')~=~f^{-1}([1,2];\{1\},\{2\}) \to 
([1,2];\{1\},\{2\})
\endaligned$$
using choices of $CW$ structures for $N,N',N''$ and crossing chain 
approximations
$$h'~:~C(\widetilde{N}') \to C(\widetilde{M}')~,~
h''~:~C(\widetilde{N}'') \to C(\widetilde{M}'')~.
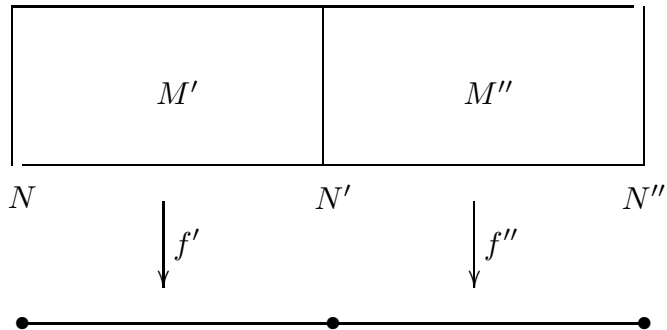
$$
$$\beginpicture
\setcoordinatesystem units <7pt,10pt>  
\putrule from -7 0 to 21 0
\putrule from -7 0 to -7 -5
\putrule from 7 0 to 7 -5
\putrule from 7 0 to 21 -5
\putrule from 21 -5 to -7 -5
\putrule from 21 0 to 21 -5
\put {$M'$} at 0 -2.7
\put {$N$} at -7 -6
\put {$N'$} at 7 -6 
\put {$M''$} at 14 -2.7
\put {$N''$} at 21 -6
\putrule from 21 -10 to -7 -10
\put {0} at -7 -11
\put {1} at 7 -11
\put {2} at 21 -11
\put {$\bullet$} at -7 -10
\put {$\bullet$} at 7 -10
\put {$\bullet$} at 21 -10
\put {$\xymatrix@R+3pt {\ar[d]^{\displaystyle{f'}} & \\ &}$} at 1.5 -7.5
\put {$\xymatrix@R+3pt {\ar[d]^{\displaystyle{f''}} & \\ &}$} at 15.5 -7.5
\endpicture$$

\vskip6pt

\noindent 
If $f'$ (resp. $f''$) has $c_i(f')$ (resp. $c_i(f'')$) critical points of
index $i$ then $f$ has 
$$c_i(f)~=~c_i(f')+c_i(f'')$$
critical points of index $i$. 

\noindent{\bf Terminology 2.5}  
Write the various chain complexes, birth, death, survival and crossing 
chain approximations for $(f',v')$ and $(f'',v'')$ as
$$\aligned
&D~=~C(\widetilde{N})~~,~~D'~=~C(\widetilde{N}')~~,~~D''~=~C(\widetilde{N}'')~,\cr
&E'~=~C(\widetilde{M}')~~,~~E''~=~C(\widetilde{M}'')~,\cr
&F'~=~C^{MS}(M',f',v')~~,~~F''~=~C^{MS}(M'',f'',v'')~,\cr
&d_{E'}~=~\pmatrix d_D & c' \cr 0 & d_{F'} \endpmatrix~:~
E'_i~=~D_i \oplus F'_i \to E_{i-1}~=~D_{i-1} \oplus F'_{i-1}~,\cr 
&d_{E''}~=~\pmatrix d_{D'} & c'' \cr 0 & d_{F''}\endpmatrix~:~
E''_i~=~D'_i \oplus F''_i \to E''_{i-1}~=~D'_{i-1} \oplus F''_{i-1}~,\cr
&h'~=~\pmatrix h'_D \cr h'_{F'} \endpmatrix~:~
D'_i \to E'_i~=~D_i \oplus F'_i~,\cr
&h''~=~\pmatrix h''_{D'} \cr h''_{F''} \endpmatrix~:~
D''_i \to E''_i~=~D'_i \oplus F''_i~.
\endaligned$$
Define the chain complexes $F,E$ by
$$\aligned
&d_F~=~\pmatrix d_{F'} & h'_{F'}c'' \cr 0 & d_{F''} \endpmatrix~:~
F_i~=~F'_i \oplus F''_i \to F_{i-1}~=~F'_{i-1} \oplus F''_{i-1}~,\cr
&d_E~=~\pmatrix d_D & c' & h'_Dc'' \cr
0 & d_{F'} & h'_{F'}c'' \cr
0 & 0 & d_{F''}\endpmatrix~:\cr
&\hskip50pt
E_i~=~D_i \oplus F'_i \oplus F''_i \to 
E_{i-1}~=~D_{i-1}\oplus F'_{i-1} \oplus F''_{i-1}~,
\endaligned$$
and let 
$$c~:~F \to D_{*-1}~~,~~h~=~\pmatrix h_D \cr h_F\endpmatrix~:~D'' \to E$$
be the chain maps defined by
$$\aligned
&c~=~\pmatrix c' & h'_Dc'' \endpmatrix~:~F_i~=~F'_i\oplus F''_i 
\to D_{i-1}~,\cr
&h_D~=~h'_Dh''_{D'}~:~D''_i \to D_i~,\cr
&h_F~=~\pmatrix h'_{F'}h''_{D'} \cr h''_{F''} \endpmatrix~:~
D''_i \to F_i~=~F'_i \oplus F''_i~.\hfill\qed
\endaligned$$

\proclaim{Proposition 2.6} 
{\rm (i)} The Morse-Smale complex of 
$$f~=~f' \cup f''~:~(M;N,N'')~=~(M';N,N') \cup (M'';N',N'')
\to ([0,2];\{0\},\{2\})$$
with respect to $v=v'\cup v'' \in \GT(f)$ is the algebraic mapping cone
$$C^{MS}(M,f,v)~=~{\Cal C}(a:F''_{*+1} \to F')$$
of the chain map $a:F''_{*+1} \to F'$ defined by
$$a~:~F''_{i+1}~=~\Z[\pi_1(M)]^{c_{i+1}(f'')} \to F'_i~=~\Z[\pi_1(M)]^{c_i(f')}~;~
\widetilde{p} \mapsto \sum\limits_{\widetilde{q}}\sum\limits_{u\in \pi_1(M)}
n(\widetilde{p},u\widetilde{q})u\widetilde{q}$$
with $n(\widetilde{p},u\widetilde{q}) \in \Z$ the algebraic number of
$\widetilde{v}$-gradient flow lines in $\widetilde{M}$ joining
$\widetilde{p}\in \widetilde{M}''$ to $u \widetilde{q}\in
\widetilde{M}'$.\newline
{\rm (ii)} For any death chain map $h'_{F'}:D' \to F'$ of $(f',v')$ and
any birth chain map $c'':F''_{*+1} \to D'$ of $(f'',v'')$ there exists 
a chain homotopy 
$$b~:~a~\simeq~h'_{F'} c''~:~F''_{*+1}  \to F'$$
such that
\roster
\item"(a)"  $b$ determines a simple isomorphism
$$I~=~\pmatrix 1 & \pm b \cr 0 & 1 \endpmatrix~:~
{\Cal C}(a)~=~C^{MS}(M,f,v) \to F~=~{\Cal C}(h'_{F'}c'')~,$$
\item"(b)" 
$cI:C^{MS}(M,f,v)_{*+1} \to D=C(\widetilde{N})$
is a birth chain map of $(f,v)$,
\item"(c)" the cellular chain complex of $\widetilde{M}$ is
$$C(\widetilde{M})~=~{\Cal C}(cI)$$
and there is defined a simple isomorphism
$$1 \oplus I~:~E \to C(\widetilde{M})~,$$
\item"(d)" 
$(1\oplus I)h:D''=C(\widetilde{N}) \to C(\widetilde{M})$
is a crossing chain approximation for $(f,v)$, with components a death chain map
$$\aligned
((1\oplus I)h)_F~=~
Ih_F~&=~\pmatrix h'_{F'}h''_{D'}\pm bh''_{F''} \cr h''_{F''}\endpmatrix~:\cr
&D''_i  \to C^{MS}(M,f,v)_i~=~F'_i \oplus F''_i
\endaligned
$$
and a survival chain map
$$((1\oplus I)h'')_D~=~h_D~=~h'_Dh''_{D'}~:~D''_i  \to D_i~.$$
\endroster
\endproclaim
\demo{Proof} (i) This is a direct consequence of the construction of the
Morse-Smale complex (2.1).\newline
(ii) Use the handlebody decomposition of $M$ (resp. $M'$, $M''$)
determined by $(f,v)$ (resp. $(f',v')$, $(f'',v'')$)
to extend the $CW$ structure on $N$ (resp. $N$, $N'$) 
to a $CW$ structure on $M$ (resp. $M'$, $M''$). 
The existence of a chain homotopy $b:a \simeq h_Fc'$ is
immediate from the observation that $a$ and $h'_{F'}c''$ are both 
connecting chain maps
for the triad of $CW$ subcomplexes $M \supset M' \supset N$ 
$$\partial ~:~F''_{*+1}~=~C(\widetilde{M},\widetilde{M}')_{*+1}
\to C(\widetilde{M}',\widetilde{N})~=~F'~,$$
which is unique up to chain homotopy.  Property (a) is just an
application of 1.2. The composite of cellular approximations
to the inclusions $N'' \to M''$, $M'' \to M$ is a cellular approximation
to the inclusion $N'' \to M$, giving (b),(c) and (d). \hfill\qed
\enddemo

\noindent{\bf Remark 2.7} 
Cornea and Ranicki \cite{CR} obtain a sharper version of 2.6 :
for any Morse map 
$$f~=~f' \cup f''~:~(M;N,N'')~=~(M;N,N') \cup (M';N',N'') \to ([0,2];\{0\},\{2\})$$
and $v=v'\cup v'' \in \GT(f)$ there exist Morse maps 
$$\widehat{f}~:~M \to \R~,~g~:~N'~=~f^{-1}(1) \to \R$$ 
and $\widehat{v} \in \GT(\widehat{f})$, $w\in \GT(g)$ such that
\roster
\item"(a)" $(\widehat{f},\widehat{v})$  agrees with $(f,v)$ outside a tubular neighbourhood of $N'$
$$f^{-1}[1-\epsilon,1+\epsilon]~=~N' \times [1-\epsilon,1+\epsilon] \subset M$$
for some small $\epsilon>0$.
\item"(b)" $(\widehat{f},\widehat{v})$ restricts to translates of $(g,w)$
$$\aligned
&(\widehat{f},\widehat{v})\vert~=~(g_+,w_+)~:~N' \times \{1+\epsilon/2\} \to \R~,\cr
&(\widehat{f},\widehat{v})\vert~=~(g_-,w_-)~:~N' \times \{1-\epsilon/2\} \to \R
\endaligned$$
with 
$$\aligned
&\hbox{\rm Crit}_i(g_+)~=~\hbox{\rm Crit}_i(g) \times \{1+\epsilon/2\}~,\cr
&\hbox{\rm Crit}_i(g_-)~=~\hbox{\rm Crit}_i(g) \times \{-\epsilon/2\}~,\cr
&\hbox{\rm Crit}_i(\widehat{f})~=~\hbox{\rm Crit}_{i-1}(g_+) 
\cup \hbox{\rm Crit}_i(g_-) \cup \hbox{\rm Crit}_i(f)~.
\endaligned$$
\item"(c)" The $v$-gradient flow lines are in one-one correspondence 
with the broken $\widehat{v}$-gradient flow lines 
i.e. joined up sequences of $\widehat{v}$-gradient flow lines which start at
critical points of $f''$ and terminate at critical points of $f'$. 
\item"(d)" The Morse-Smale complex of $(\widehat{f},\widehat{v})$ 
is of the form
$$\aligned
&d_{C^{MS}(M,\widehat{f},\widehat{v})}~=~
\pmatrix d_{F'} & h'_{F'} & 0 & 0 \cr
0 & -d_{D'} & 0 & 0 \cr
0 & 1 & d_{D'} & c'' \cr
0 & 0 & 0 & d_{F''} \endpmatrix~:\cr
&C^{MS}_i(M,\widehat{f},\widehat{v})~=~F'_i \oplus D'_{i-1} \oplus D'_i
\oplus F''_i \cr
&\hskip100pt
\to C^{MS}_{i-1}(M,\widehat{f},\widehat{v})~=~F'_{i-1} \oplus D'_{i-2} \oplus D'_{i-1}
\oplus F''_{i-1} 
\endaligned$$
with $D'=C(N',g,w)$, giving choices of 'gradient-like'
crossing chain approximations $h',h''$  such that the chain homotopy in 2.6 (ii) 
is 
$$\aligned
&b~=~0~:~a~ \simeq~ h'_{F'}c''~:\cr
&F''~=~C^{MS}(M'',f'',v'') \to F'_{*-1}~=~C^{MS}(M',f',v')_{*-1}
\endaligned$$
(i.e. $a=h'_{F'}c''$) and the simple isomorphism
$$I~:~C^{MS}(M,f,v)~=~{\Cal C}(a) \to F~=~{\Cal C}(h'_{F'}c'')$$
is the identity.\hfill\qed
\endroster

\subheading{3. The proper Morse-Smale complex}

The construction of $C^{MS}(M,f,v)$ applies just as well to a 
proper real-valued Morse function on a non-compact manifold :

\noindent{\bf Definition 3.1} 
Let $(M,\partial M)$ be a non-compact manifold with compact boundary and 
a proper real-valued Morse function
$$f~:~(M,\partial M) \to ([0,\infty),\{0\})~,$$
and let $v \in \GT(f)$. The
{\it proper Morse-Smale complex} $C^{MS}(M,f,v)$ is defined 
exactly as in the compact case, with
\roster
\item"{\rm (i)}" $C^{MS}_i(M,f,v)$ the based free $\Z[\pi_1(M)]$-module generated
by $\gamma_i(f)$, the set of critical points of $f$ with index $i$,
\item"{\rm (ii)}" the boundary $\Z[\pi_1(M)]$-module morphisms are given by
$$\aligned
d~:~C^{MS}_i(M,f,v)~=~\Z[\pi_1(M)]^{\gamma_i(f)}& \to 
C^{MS}_{i-1}(M,f,v)~=~\Z[\pi_1(M)]^{\gamma_{i-1}(f)}~;\cr
&\widetilde{p} \mapsto \sum\limits_{\widetilde{q}}\sum\limits_{u\in \pi_1(M)}
n(\widetilde{p},u\widetilde{q})u\widetilde{q}
\endaligned$$
with $n(\widetilde{p},u\widetilde{q}) \in \Z$ the algebraic number of 
$v$-gradient flow lines in $\widetilde{M}$ joining $\widetilde{p}$ to 
$u \widetilde{q}$.
\endroster
$$\eqno{\qed}$$

Given a proper Morse function 
$f:(M,\partial M) \to ([0,\infty),\{0\})$ and a gradient-like vector field
for $f$ and a $CW$ structure for $\partial M$ let $M$ have the $CW$ structure
given by the handle decomposition
$$M~=~\partial M \times I \cup\bigcup^m_{i=0}\bigcup_{\gamma_i(f)}h^i~.$$
The expression of 2.6 for the Morse-Smale complex of the union of two
Morse functions on adjoining compact cobordisms will now be applied to 
obtain an isomorphism between the Morse-Smale complex $C^{MS}(M,f,v)$ and 
the relative cellular chain complex of $(\widetilde{M},\partial\widetilde{M})$ 
$$I~:~C^{MS}(M,f,v)~\cong~C(\widetilde{M},\partial\widetilde{M})~.$$

\noindent{\bf Terminology 3.2} For $j=0,1,2,\dots$ let
$$\aligned
&f[j]~=~f\vert~:\cr
&(M[j];N[j],N[j+1])~=~f^{-1}([j,j+1];\{j\},\{j+1\})\to ([j,j+1];\{j\},\{j+1\})
\endaligned$$
be the Morse functions on compact cobordisms given by the restriction of $f$.
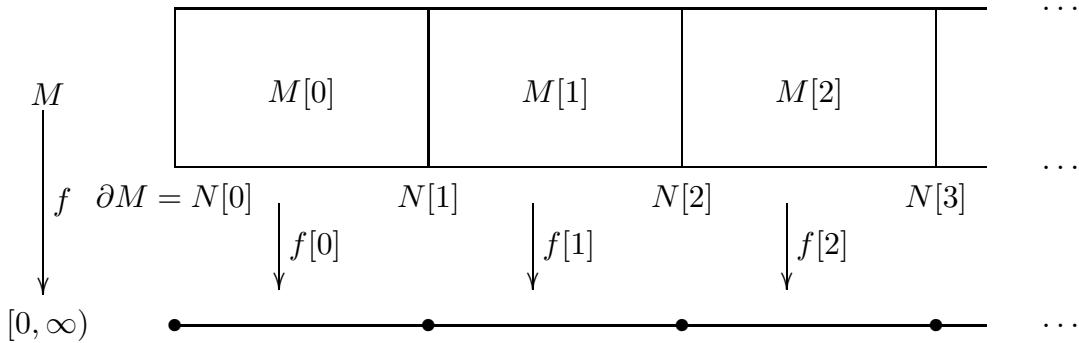
$$\beginpicture
\setcoordinatesystem units <8pt,10pt>  
\putrule from -15 0 to 17 0
\putrule from -15 -5 to 17 -5
\putrule from -15 0 to -15 -5
\putrule from -5 0 to -5 -5
\putrule from 5 0 to 5 -5
\putrule from 15 0 to 15 -5
\putrule from -15 -10 to 17 -10
\put {$M[0]$} at -10 -2.7
\put {$M[1]$} at 0 -2.7
\put {$M[2]$} at 10 -2.7 
\put {$\partial M=N[0]$} at -15 -6
\put {$N[1]$} at -5 -6
\put {$N[2]$} at 5 -6 
\put {$N[3]$} at 15 -6
\put {$0$} at -15 -11
\put {$1$} at -5 -11
\put {$2$} at 5 -11
\put {$3$} at 15 -11
\put {$\bullet$} at -15 -10
\put {$\bullet$} at -5 -10
\put {$\bullet$} at 5 -10
\put {$\bullet$} at 15 -10
\put {$M$} at -20 -2.7
\put {$[0,\infty)$} at -20 -10
\put {$\xymatrix@R+34pt {\ar[d]^{\displaystyle f}&\\&}$} at -18.3 -6.1
\put {$\xymatrix@R+3pt {\ar[d]^{\displaystyle{f[0]}}&\\&}$} at -9 -7.5
\put {$\xymatrix@R+3pt {\ar[d]^{\displaystyle{f[1]}}&\\&}$} at 1 -7.5
\put {$\xymatrix@R+3pt {\ar[d]^{\displaystyle{f[2]}}&\\&}$} at 11 -7.5
\put {$\dots$} at  20 0
\put {$\dots$} at  20 -5
\put {$\dots$} at  20 -10
\endpicture$$
\vskip4pt

\noindent The inclusions
$$g[j]~:~N[j] \to M[j]~,~h[j]~:~N[j+1] \to M[j]$$
induce embeddings of based f.\,g. free $\Z[\pi_1(M)]$-module chain complexes
$$g[j]~:~C(\widetilde{N}[j]) \to C(\widetilde{M}[j])$$
and chain maps
$$h[j]~:~C(\widetilde{N}[j+1]) \to C(\widetilde{M}[j])~.$$
Write
$$\aligned
&D[j]~=~C(\widetilde{N}[j])~~,~~E[j]~=~C(\widetilde{M}[j])~~,~~
F[j]~=~C(\widetilde{M}[j],\widetilde{N}[j])~~,\cr
&d_{E[j]}~=~\pmatrix d_{D[j]} & c[j] \cr 0 & d_{F[j]} \endpmatrix~:\cr
&\hskip15mm
E_i[j]~=~D_i[j]\oplus F_i[j] \to E_{i-1}[j]~=~D_{i-1}[j]\oplus F_{i-1}[j]~,\cr
&g[j]~=~\pmatrix 1 \cr 0 \endpmatrix~ :~ D_i[j] \to E_i[j]~=~D_i[j]\oplus F_i[j]~,\cr
&h[j+1]~=~\pmatrix h_{D[j+1]} \cr h_{F[j+1]} \endpmatrix : D_i[j+1] \to E_i[j]~=~D_i[j]\oplus F_i[j]~.
\hfill \qed
\endaligned$$

\proclaim{\bf Proposition 3.3} Let $f:(M,\partial M) \to ([0,\infty),\{0\})$
be a proper Morse function, and let $v \in \GT(f)$.\newline
{\rm (i)} The cellular chain complex of $(M,\partial M)$ is of the form
$$C(\widetilde{M},\partial\widetilde{M})~=~\hbox{\rm coker}
\big(g-h:\sum\limits^{\infty}_{j=0}C(\widetilde{N}[j])
\to \sum\limits^{\infty}_{j=0}C(\widetilde{M}[j])\big)~,$$
a based free $\Z[\pi_1(M)]$-module complex
with basis the images of the basis elements
in $\sum\limits^{\infty}_{j=0}C(\widetilde{M}[j],\widetilde{N}[j])$,
and may be expressed as 
$$\aligned
&d_{C(\widetilde{M},\partial\widetilde{M})}~=~
\pmatrix d_{F[0]} & h_{F[0]}c[1] & h_{F[0]}h_{D[1]}c[2] & 
h_{F[0]}h_{D[1]}h_{D[2]}c[3] & \dots \cr\vspace{2mm}
0 & d_{F[1]} & h_{F[1]}c[2] & h_{F[1]}h_{D[2]}c[3] & \dots \cr\vspace{2mm}
0 & 0 & d_{F[2]} & h_{F[2]}c[3] & \dots \cr\vspace{2mm}
\vdots & \vdots & \vdots & \vdots & \ddots \endpmatrix\cr
&\hskip5mm 
:~C_i(\widetilde{M},\partial\widetilde{M})~=~
\sum\limits^{\infty}_{k=0}F_i[k] \to 
C_{i-1}(\widetilde{M},\partial\widetilde{M})~=~
\sum\limits^{\infty}_{j=0}F_{i-1}[j]~.\endaligned$$
{\rm (ii)} The Morse-Smale complex of $(f,v)$ is of the form
$$\aligned
&d_{C^{MS}(M,f,v)}~=\cr
&\pmatrix d_{C^{MS}(M[0],f[0],v[0])} & c^{MS}[0,1] & c^{MS}[0,2] & 
c^{MS}[0,3] & \dots \cr\vspace{2mm}
0 & d_{C^{MS}(M[1],f[1],v[1])} & c^{MS}[1,2] & c^{MS}[1,3] & \dots \cr\vspace{2mm}
0 & 0 & d_{C^{MS}(M[2],f[2],v[2])} & c^{MS}[2,3] & \dots \cr\vspace{2mm}
\vdots & \vdots & \vdots & \vdots & \ddots \endpmatrix\cr
&:~C^{MS}_i(M,f,v)~=~\sum\limits^{\infty}_{k=0}C^{MS}_i(M[k],f[k],v[k])\cr
&\hskip120pt
\to C^{MS}_{i-1}(M,f,v)~=~\sum\limits^{\infty}_{j=0}C^{MS}_{i-1}(M[j],f[j],v[j])
\endaligned$$
with 
$$c^{MS}[j,k]~:~C_i^{MS}(M[k],f[k],v[k]) \to C^{MS}_{i-1}(M[j],f[j],v[j])~~(j<k)$$
counting the $\widetilde{v}$-gradient flow lines in the universal cover $\widetilde{M}$
which start at an index $i$ critical point of $\widetilde{f}[k]$ and
terminate at an index $i-1$ critical point of $\widetilde{f}[j]$.\newline
{\rm (iii)} There exists an isomorphism of chain complexes
$$I~:~C^{MS}(M,f,v)~\cong~C(\widetilde{M},\partial\widetilde{M})$$
of the form
$$I~=~1+ \sum\limits_{j'<j}b[j',j]~:~
C^{MS}_i(M,f,v)~=~\sum\limits^{\infty}_{j=0}F_i[j]
\to C_i(\widetilde{M},\partial\widetilde{M})~=~
\sum\limits^{\infty}_{j'=0}F_i[j']~.$$
\endproclaim
\demo{Proof} (i) The exact sequence 
$$0 \to \sum\limits^{\infty}_{j=0}C(\widetilde{N}[j]) 
\xymatrix{\ar[r]^{\displaystyle{g-h}}&}  
\sum\limits^{\infty}_{j=0}C(\widetilde{M}[j]) \to C(\widetilde{M},
\partial\widetilde{M}) \to 0$$
is just the chain level Mayer-Vietoris sequence for the union
$$M~=~M[\hbox{even}] \cup M[\hbox{odd}]$$
with
$$M[\hbox{even}]~=~\bigcup\limits_{j~\hbox{even}}M[j]~~,~~
M[\hbox{odd}]~=~\bigcup\limits_{j~\hbox{odd}}M[j]~.$$
The matrix formula for $\hbox{\rm coker}(g-h)$ is a direct application
of 1.8 (i) (b) with 
$$\phi~=~g-h~=~\pmatrix 1-h_D \cr -h_F \endpmatrix~:~D_i \to E_i~=~D_i \oplus F_i~,$$
noting that each
$$1-h_D~=~
\pmatrix 1 & -h_{D[0]} & 0 & \dots \cr
0 & 1 & -h_{D[1]} & \dots \cr
0 & 0 & 1 & \dots \cr
\vdots & \vdots & \vdots & \ddots \endpmatrix~:~
E_i~=~\sum\limits^{\infty}_{j=0}F_i[j] \to E_i~=~\sum\limits^{\infty}_{j=0}F_i[j]$$
is an automorphism, with inverse
$$\aligned
&(1-h_D)^{-1}~=~
\pmatrix 1 & h_{D[0]} & h_{D[0]}h_{D[1]} & \dots \cr\vspace{1.5ex}
0 & 1 & h_{D[1]} & \dots \cr\vspace{1.5ex}
0 & 0 & 1 & \dots \cr\vspace{1.5ex}
\vdots & \vdots & \vdots & \ddots \endpmatrix~:\cr
&\hskip25pt
E_i~=~\sum\limits^{\infty}_{j=0}F_i[j] \to E_i~=~\sum\limits^{\infty}_{j=0}F_i[j]~.
\endaligned$$
(ii) By construction.\newline
(iii) For $k=1,2,\dots$ define the Morse function
$$f[0,k]~=~\bigcup^k_{j=1}f[j]~=~f\vert~:~
M[0,k]~=~\bigcup^k_{j=1}M[j]~=~f^{-1}[0,k] \to [0,k]~,$$
and assume inductively that there is given an isomorphism
$$I[0,k]~:~C^{MS}(M[0,k],f[0,k],v[0,k])~\cong~C(\widetilde{M}[0,k],\partial\widetilde{M})$$
of the form
$$\aligned
&I[0,k]~=~1+ \sum\limits_{j'<j}b[j',j]~:\cr
&C^{MS}_i(M[0,k],f[0,k],v[0,k])~=~\sum\limits^k_{j=0}F_i[j]
\to C_i(\widetilde{M}[0,k],\partial\widetilde{M})~=~
\sum\limits^k_{j'=0}F_i[j']~.
\endaligned$$
Now apply 2.6 (i) to the Morse function
$$\aligned
&f[0,k+1]~=~f[0,k] \cup f[k]~:\cr
&M[0,k+1]~=~M[0,k] \cup M[k] \to [0,k+1]~=~[0,k] \cup [k,k+1]~,
\endaligned$$
extending $I[0,k]$ to an isomorphism $I[0,k+1]$ of the same form, and
pass to the direct limit to obtain an isomorphism
$$\aligned
&I~=~\varinjlim_k I[0,k]~:~
C^{MS}(M,f,v)~=~\varinjlim_k C^{MS}(M[0,k],f[0,k],v[0,k])\cr
& \hskip150pt \to C(\widetilde{M},\partial\widetilde{M})~=~
\varinjlim_k C(\widetilde{M}[0,k],\partial\widetilde{M})~.
\endaligned$$
of the form
$$I~=~1+\sum\limits_{j'<j}b[j',j]~.\eqno{\qed}$$
\enddemo

\subheading{4. The Cohn noncommutative localization}

We refer to Cohn \cite{C} and Schofield \cite{Scho} for general accounts 
of the localization $\Sigma^{-1}R$ of a ring $R$ inverting a
set $\Sigma$ of square matrices. The natural morphism $R \to \Sigma^{-1}R$
has the universal property that a morphism of rings
$R \to A$ which sends $\Sigma$ to invertible matrices in $A$ has a unique
factorization
$$R \to \Sigma^{-1}R \to A~.$$
The Gerasimov-Malcolmson normal form expresses every morphism 
of f.\,g. free $\Sigma^{-1}R$-modules
$$\phi~:~\Sigma^{-1}R^n \to \Sigma^{-1}R^p$$
as a composite 
$$\phi~=~f\sigma^{-1}g~:~\Sigma^{-1}R^n \to 
\Sigma^{-1}R^m \to \Sigma^{-1}R^m \to \Sigma^{-1}R^p$$
(nonuniquely) for some $R$-module morphisms 
$$f~:~R^m \to R^p~~,~~\sigma~:~R^m \to R^m~~,~~g~:~R^n \to R^m$$
such that $\sigma$ is $\Sigma^{-1}R$-invertible.

\proclaim{\bf Proposition 4.1}  {\rm (Sheiham \cite{Sh,3.1})}
Given a ring morphism $\epsilon:R \to A$ let $\Sigma$ 
be the set of all square matrices in $R$ which become invertible in $A$.\newline
{\rm (i)} The ring morphism $\epsilon$ extends to a ring morphism
$$\epsilon~:~\Sigma^{-1}R \to A~.$$
{\rm (ii)} An endomorphism of a f.\,g. free $\Sigma^{-1}R$-module
$$\phi~:~\Sigma^{-1}R^n \to \Sigma^{-1}R^n$$
is an automorphism if and only if $\epsilon(\phi):A^n \to A^n$ is an
$A$-module automorphism.
\endproclaim
\demo{Proof} 
(i) By the universal property of $R \to \Sigma^{-1}R$.\newline
(ii) It is clear that if $\phi$ is an automorphism then so is
$\epsilon(\phi)$.\newline
\indent Conversely, suppose that $\epsilon(\phi)$ is an automorphism.
Express $\phi$ in the Gerasimov-Malcolmson normal form
$$\phi~=~f\sigma^{-1}g~:~\Sigma^{-1}R^n \to \Sigma^{-1}R^n$$
for some $R$-module morphisms 
$$f~:~R^m \to R^n~~,~~\sigma~:~R^m \to R^m~~,~~g~:~R^n \to R^m$$
such that $\epsilon(\sigma):A^m \to A^m$ is an automorphism.
The $R$-module endomorphism defined by
$$\theta~=~\pmatrix 0 & -f \cr g & \sigma \endpmatrix~:~
R^n \oplus R^m \to R^n \oplus R^m$$
is $A$-invertible, since the induced $\Sigma^{-1}R$-module endomorphism
$$\theta~=~\pmatrix \phi & -f \cr 0 & \sigma \endpmatrix
\pmatrix 1 & 0 \cr \sigma^{-1}g & 1 \endpmatrix~:~
\Sigma^{-1}(R^n \oplus R^m) \to \Sigma^{-1}(R^n \oplus R^m)$$
is $A$-invertible. Thus $\theta:\Sigma^{-1}(R^n\oplus R^m) \to
\Sigma^{-1}(R^n\oplus R^m)$ is a $\Sigma^{-1}R$-module automorphism,
and hence so is $\phi$.
\hfill\qed
\enddemo

In the application of Cohn localization in \S5 below 
$$\epsilon~:~R~=~A_{\alpha }[z] \to A~;~z \mapsto 0~.$$

\subheading{5. Polynomial extensions}

The relationship between the algebraic mapping cone and cokernel worked
out in \S1 will now be applied to the algebraic situation arising from
a circle-valued Morse function $f:M \to S^1$.  The actual application
to the Novikov complex $C^{Nov}(M,f,v)$ will be carried out in \S6.

\noindent{\bf Definition 5.1}
Let $A$ be a ring with an automorphism $\alpha :A \to A$, and let
$z$ be an indeterminate over $A$ such that 
$$az~=~z\alpha (a)~~(a \in A)~.$$
{\rm (i)} The {\it $\alpha $-twisted Laurent polynomial extension} of $A$ 
$$A_{\alpha }[z,z^{-1}]~=~\sum\limits^{\infty}_{j=-\infty}z^jA~=~
A_{\alpha }[z,z^{-1}]$$
is the ring of polynomials 
$\sum\limits^{\infty}_{j=-\infty}a_jz^j$ $(a_j \in A)$
such that $\{j \in \Z\,\vert\,a_j \neq 0\}$ is finite.\newline
{\rm (ii)} The {\it $\alpha $-twisted Novikov completion of $A_{\alpha }[z,z^{-1}]$}
$$A_{\alpha }((z))~
=~\varprojlim_k \sum\limits^k_{j=-\infty}z^jA~=~A_{\alpha }[[z]][z^{-1}]$$
is the ring of formal power series $\sum\limits^{\infty}_{j=-\infty}a_jz^j$ $(a_j \in A)$
such that $\{j \leq 0\,\vert\,a_j \neq 0\}$ is finite.\hfill\qed

Given an $A$-module $B$ and $j \in \Z$ let $z^jB$ be the $A$-module
with elements $z^jx$ ($x \in B$) and 
$$a(z^jx)~=~z^j\alpha ^j(a)x~~,~~z^jx+z^jx'~=~z^j(x+x')~~(a\in A, x,x' \in B)~.$$
The induced $A_{\alpha }[z,z^{-1}]$-module is then given by 
$$A_{\alpha }[z,z^{-1}]\otimes_AB~=~B_{\alpha }[z,z^{-1}]~=~\sum\limits^{\infty}_{j=-\infty}z^jB~.$$
For any $A$-module $C$ and $k \in \Z$ the $A$-module morphisms 
$z^j B \to z^kC$ are given by
$$z^{k-j}\theta ~:~ z^jB \to z^kC ~;~z^jx \mapsto z^k\theta(x)$$
with $\theta:B \to C$ a morphism of the additive groups such that
$$\theta(ax)~=~\alpha ^{k-j}(a)\theta(x) \in C ~~(a \in A, x \in B)~.$$

We shall write $z^{-1}B$ as $\alpha B$.

For a f.\,g. free $A$-module $B$ and any $A$-module $C$
every $A_{\alpha }[z,z^{-1}]$-module morphism
$$\psi~:~B_{\alpha }[z,z^{-1}] \to C_{\alpha }[z,z^{-1}]$$
is given by
$$\psi~=~\sum\limits^{\infty}_{j=-\infty}z^j\psi_j~:~
\sum\limits^{\infty}_{k=-\infty}z^kx_k \mapsto 
\sum\limits^{\infty}_{j=-\infty}\sum\limits^{\infty}_{k=-\infty}z^{j+k}
\psi_j(x_k)$$
with $\psi_j:B \to z^jC$ $A$-module morphisms such that $\{j \in \Z\,\vert\,
\psi_j \neq 0\}$ is finite. Similarly, 
every $A_{\alpha }((z))$-module morphism
$$\widehat{\psi}~:~B_{\alpha }((z)) \to C_{\alpha }((z))$$
is given by
$$\widehat{\psi}~=~\sum\limits^{\infty}_{j=-\infty}z^j\widehat{\psi}_j~:~
\sum\limits^{\infty}_{k=-\infty}z^kx_k \mapsto 
\sum\limits^{\infty}_{j=-\infty}\sum\limits^{\infty}_{k=-\infty}z^{j+k}
\widehat{\psi}_j(x_k)$$
with $\widehat{\psi}_j:B \to z^jC$ $A$-module morphisms such that 
$\{j \leq 0\,\vert\, \widehat{\psi}_j \neq 0\}$ is finite.

\noindent{\bf Definition 5.2} Let $\Sigma^{-1}A_{\alpha }[z,z^{-1}]$ be
the localization of $A_{\alpha }[z,z^{-1}]$ inverting the set $\Sigma$
of square matrices in $A_{\alpha }[z] \subset A_{\alpha }[z,z^{-1}]$
which become invertible in $A$ under the augmentation 
$$\epsilon~:~A_{\alpha }[z] \to A~;~z \mapsto 0~.\eqno{\qed}$$
$$\beginpicture
\setcoordinatesystem units <8pt,10pt>  
\putrule from -17 0 to 17 0
\putrule from -17 -5 to 17 -5
\putrule from -17 0 to -15 -5
\putrule from -15 0 to -15 -5
\putrule from -5 0 to -5 -5
\putrule from 5 0 to 5 -5
\putrule from 15 0 to 15 -5
\put {$zE$} at -10 -2.7
\put {$E$} at 0 -2.7
\put {$z^{-1}E$} at 10 -2.7 
\put {$zD$} at -15 -6
\put {$D$} at -5 -6
\put {$z^{-1}D=\alpha D$} at 5 -6 
\put {$z^{-2}D$} at 15 -6
\put {$\dots$} at  20 0
\put {$\dots$} at  20 -5
\put {$\dots$} at  -20 0
\put {$\dots$} at  -20 -5
\endpicture$$

\proclaim{\bf Proposition 5.3} 
Let $D,E$ be f.\,g. free $A$-module chain complexes, and let 
$g:D \to E$, $h:\alpha D \to E$ be $A$-module chain maps such that
$$\aligned
&g~=~\pmatrix 1 \cr 0 \endpmatrix~:~D_i \to E_i ~=~D_i \oplus F_i~,\cr
&h~=~\pmatrix h_D \cr h_F \endpmatrix~:~\alpha D_i \to E_i ~=~D_i \oplus F_i~,\cr
&d_E~=~\pmatrix d_D & c \cr 0 & d_F \endpmatrix~:~
E_i~=~D_i \oplus F_i \to E_{i-1}~=~D_{i-1} \oplus F_{i-1}~.
\endaligned$$
Given bases for $D,F$ let $E$ have the corresponding basis.\newline
{\rm (i)} The $\Sigma^{-1}D_{\alpha }[z,z^{-1}]$-module morphism
$$\phi~=~g-zh~:~\Sigma^{-1}D_{\alpha }[z,z^{-1}] \to \Sigma^{-1}E_{\alpha }[z,z^{-1}]$$
is an embedding, the natural projection
$$p~:~C~=~{\Cal C}(\phi) \to \widehat{C}~=~\hbox{\rm coker}(\phi)$$
is a chain equivalence, and the inclusions $F_i \to E_i$ induce an isomorphism
$$\widehat{F}~\cong~\widehat{C}$$ 
with $\widehat{F}$ the $\Sigma^{-1}A_{\alpha }[z,z^{-1}]$-module chain complex 
given by
$$\aligned
d_{\widehat{F}}~=~~&d_F+zh_F(1-zh_D)^{-1}c~:\cr
& \widehat{F}_i~=~\Sigma^{-1}(F_i)_{\alpha }[z,z^{-1}] 
\to\widehat{F}_{i-1}~=~\Sigma^{-1}(F_{i-1})_{\alpha }[z,z^{-1}]~.
\endaligned$$
If $\widehat{C}$ is given the basis determined by the bases of $D,F$ 
and the isomorphism $\widehat{F}\cong\widehat{C}$ then
$$\aligned
\tau(p:C \simeq \widehat{C})~=~&-\sum\limits^{\infty}_{i=0}(-)^i
\tau(1-zh_D:
\Sigma^{-1}(D_i)_{\alpha }[z,z^{-1}] \to\Sigma^{-1}(D_i)_{\alpha }[z,z^{-1}])\cr
& \in K_1( \Sigma^{-1}A_{\alpha }[z,z^{-1}] )~.
\endaligned$$
{\rm (ii)} An $A$-module chain homotopy
$$k~:~h~ \simeq~ h'~:~\alpha D \to E$$
determines a $\Sigma^{-1}A_{\alpha }[z,z^{-1}]$-module chain isotopy 
$$\psi~:~\phi~=~g-zh~\sim~\phi'~=~g-zh'~:~\Sigma^{-1}D_{\alpha }[z,z^{-1}] \to 
\Sigma^{-1}E_{\alpha }[z,z^{-1}]$$
and simple isomorphisms of $\Sigma^{-1}A_{\alpha }[z,z^{-1}]$-module chain complexes
$$C~\cong~C'~~,~~\widehat{F}~\cong~\widehat{F}'~~,~~\widehat{C}~\cong~\widehat{C}'$$
where $C={\Cal C}(\phi)$, $C'={\Cal C}(\phi')$ etc., with
$$\aligned
&\tau(C\cong C')~=~0 \in K_1(\Sigma^{-1}A_{\alpha }[z,z^{-1}])~,\cr
&\tau(\widehat{C}\cong \widehat{C}')~=~
\sum\limits^{\infty}_{i=0}(-)^i\tau(1-zh_D:
\Sigma^{-1}(D_i)_{\alpha }[z,z^{-1}] \to \Sigma^{-1}(D_i)_{\alpha }[z,z^{-1}])\\
&\hphantom{\tau(\widehat{C}\cong \widehat{C}')~=~}
-\sum\limits^{\infty}_{i=0}(-)^i\tau(1-zh'_D:
\Sigma^{-1}(D_i)_{\alpha }[z,z^{-1}] \to \Sigma^{-1}(D_i)_{\alpha }[z,z^{-1}])\cr
&\hphantom{\tau(\widehat{C}\cong \widehat{C}')~}
=~0 \in K_1(\Sigma^{-1}A_{\alpha }[z,z^{-1}])~.
\endaligned$$
{\rm (iii)} Let $\{D(k)\}$, $\{E(k)\}$, $\{F(k)\}$ be the inverse systems of 
$A$-module chain complexes defined by
$$\aligned
&D(k)~=~D_{\alpha }[z,z^{-1}]/\sum\limits^{\infty}_{j=k+1}z^jD~
=~\sum\limits^k_{j=-\infty}z^jD~,\cr
&E(k)~=~E_{\alpha }[z,z^{-1}]/\sum\limits^{\infty}_{j=k+1}z^jE~
=~\sum\limits^k_{j=-\infty}z^jE~,\cr
&F(k)~=~\hbox{\rm coker}\big(g-zh:D(k) \to E(k)\big)
\endaligned$$
with structure maps the natural projections
$$D(k)\to D(k-1)~~,~~E(k)\to E(k-1)~~,~~F(k)\to F(k-1)~.$$
The short exact sequence of inverse systems
of $A$-module chain complexes
$$0 \to D(k) \raise4pt\hbox{$g-zh \atop\to$} E(k) \to F(k) \to 0$$
induces a short exact sequence of the inverse limit
$A_{\alpha }((z))$-module chain complexes
$$0 \to \varprojlim_kD(k) \raise4pt\hbox{$g-zh \atop\to$} 
\varprojlim_kE(k) \to \varprojlim_kF(k) \to 0$$
with 
$$\varprojlim_k D(k)~=~D_{\alpha }((z))~~,~~\varprojlim_k E(k)~=~E_{\alpha }((z))~.$$
Moreover, the inclusions $F_i \to E_i$ induce isomorphisms
$$A_{\alpha }((z))\otimes_{\Sigma^{-1}A_{\alpha }[z,z^{-1}]}\widehat{F}~
\cong~\varprojlim_k F(k)~.$$
\endproclaim
\demo{Proof} 
(i) Immediate from 1.8 (i), since
$$\aligned
&\phi~=~\pmatrix \phi_D \cr \phi_F \endpmatrix~=~
\pmatrix 1-zh_D \cr -zh_F \endpmatrix~:\cr
&\Sigma^{-1}(D_i)_{\alpha }[z,z^{-1}] \to \Sigma^{-1}(E_i)_{\alpha }[z,z^{-1}]~=~
\Sigma^{-1}(D_i\oplus F_i)_{\alpha }[z,z^{-1}]
\endaligned$$
with
$$\phi_D~=~1-zh_D~:~\Sigma^{-1}(D_i)_{\alpha }[z,z^{-1}] \to \Sigma^{-1}(D_i)_{\alpha }[z,z^{-1}]$$
a $\Sigma^{-1}A_{\alpha }[z,z^{-1}]$-module automorphism.\newline
(ii) Write the $A$-module chain homotopy $k:h\simeq h':\alpha D \to E$ as
$$k~=~\pmatrix k_D \cr k_F \endpmatrix~:~
\alpha D_i \to E_{i+1}~=~D_{i+1}\oplus F_{i+1}~,$$
so that
$$\pmatrix h'_D \cr h'_F \endpmatrix - \pmatrix h_D \cr h_F \endpmatrix ~=~
\pmatrix d_D & c \cr 0 & d_F \endpmatrix 
\pmatrix k_D \cr k_F \endpmatrix + \pmatrix k_D \cr k_F \endpmatrix d_D~:~
\alpha D_i \to E_i~=~D_i \oplus F_i~.$$
Define a $\Sigma^{-1}A_{\alpha }[z,z^{-1}]$ chain homotopy
$$\theta~:~\phi~\simeq~\phi'~:~\Sigma^{-1}D_{\alpha }[z,z^{-1}]
\to \Sigma^{-1}E_{\alpha }[z,z^{-1}]$$
by
$$\aligned
&\theta~=~\pmatrix \theta_D \cr \theta_F \endpmatrix~=~
\pmatrix -zk_D \cr -zk_F \endpmatrix~:\cr
&\hskip20pt
\Sigma^{-1}(D_i)_{\alpha }[z,z^{-1}]
\to \Sigma^{-1}(E_{i+1})_{\alpha }[z,z^{-1}]~=~
\Sigma^{-1}(D_{i+1}\oplus F_{i+1})_{\alpha }[z,z^{-1}]~.
\endaligned$$
As in 1.8 (ii) the $\Sigma^{-1}A_{\alpha }[z]$-module morphisms
$$\aligned
&\psi~=~\pmatrix \theta_D(\phi_D)^{-1} & 0 \cr 
\theta_F(\phi_D)^{-1} & 0 \endpmatrix~=~
\pmatrix -zk_D(1-zh_D)^{-1} & 0 \cr 
-zk_F(1-zh_D)^{-1} & 0 \endpmatrix~:\cr
&\hskip20pt
\Sigma^{-1}(E_i)_{\alpha }[z]~=~\Sigma^{-1}(D_i\oplus F_i)_{\alpha }[z]
\to \Sigma^{-1}(E_i)_{\alpha }[z]~=~\Sigma^{-1}(D_i\oplus F_i)_{\alpha }[z]
\endaligned$$
are such that
$$\phi'~=~(1+d_E\psi+\psi d_E)\phi~:~D \to E~.$$
The $\Sigma^{-1}A_{\alpha }[z]$-module endomorphism
$$\aligned
&1+d_E\psi+\psi d_E~=~
\pmatrix 1 & 0 \cr 0 & 1 \endpmatrix
+ \pmatrix d_D & c \cr 0 & d_F \endpmatrix
\pmatrix -zk_D(1-zh_D)^{-1} & 0 \cr 
-zk_F(1-zh_D)^{-1} & 0 \endpmatrix\cr
&\hphantom{1+d_E\psi+\psi d_E~=~\pmatrix 1 & 0 \cr 0 & 1 \endpmatrix}  
+ \pmatrix -zk_D(1-zh_D)^{-1} & 0 \cr 
-zk_F(1-zh_D)^{-1} & 0 \endpmatrix\pmatrix d_D & c \cr 0 & d_F \endpmatrix\cr
&\hskip100pt
:~\Sigma^{-1}(E_i)_{\alpha }[z] \to \Sigma^{-1}(E_i)_{\alpha }[z]
\endaligned$$
has augmentation an $A$-module automorphism
$$\epsilon(1+d_E\psi+\psi d_E)~=~1~:~E_i \to E_i~,$$
so that 
$$1+d_E\psi+\psi d_E~:~\Sigma^{-1}(E_i)_{\alpha }[z]\to
\Sigma^{-1}(E_i)_{\alpha }[z]$$
is a $\Sigma^{-1}A_{\alpha }[z]$-module automorphism by 4.1,
and $\psi$ defines a chain isotopy
$$\psi~:~\phi~\sim~\phi'~:~
\Sigma^{-1}D_{\alpha }[z] \to \Sigma^{-1}E_{\alpha }[z]~.$$
Define the isomorphisms
$$\aligned
&q~=~\pmatrix 1 & \pm \psi\phi  \cr 0 & 1 \endpmatrix~:~
C~=~{\Cal C}(\phi) \to C'~=~{\Cal C}(\phi')~,\cr
&r~=~[1+d_E\psi+\psi d_E]~:~\widehat{C}~=~\hbox{\rm coker}(\phi) \to 
\widehat{C}'~=~\hbox{\rm coker}(\phi')
\endaligned$$
as in 1.5, with $q$ simple. By 1.8 (ii) and (i) 
$$\aligned
\tau(r)~&=~\tau(p:C \simeq \widehat{C}) -\tau(p':C' \simeq \widehat{C}')\cr
&=~\sum\limits^{\infty}_{i=0}(-)^i\tau(1-zh_D:
\Sigma^{-1}(D_i)_{\alpha }[z,z^{-1}] \to \Sigma^{-1}(D_i)_{\alpha }[z,z^{-1}])\\
&\hphantom{\tau(\widehat{C}\cong \widehat{C}')~=~}
-\sum\limits^{\infty}_{i=0}(-)^i\tau(1-zh'_D:
\Sigma^{-1}(D_i)_{\alpha }[z,z^{-1}] \to \Sigma^{-1}(D_i)_{\alpha }[z,z^{-1}])\cr
&\hphantom{\tau(\widehat{C}\cong \widehat{C}')~}
=~0 \in K_1(\Sigma^{-1}A_{\alpha }[z,z^{-1}])~.
\endaligned$$
(iii) The $A$-module chain map
$$\phi(k)~=~g-zh~:~D(k) \to E(k)$$
is of the type considered in 1.8, with the components of
$$\phi(k)~=~\pmatrix \phi_D(k) \cr \phi_F(k) \endpmatrix~:~D(k)_i \to E(k)_i~=~D(k)_i \oplus F(k)_i$$
given by
$$\aligned
&\phi_D(k)~=~\pmatrix 1 & -zh_D & 0 & \dots \cr 
               0 & 1 & -zh_D & \dots \cr 
               0 & 0 & 1 & \dots \cr
               \vdots & \vdots & \vdots & \ddots \endpmatrix~:\cr
&\hskip100pt
D(k)_i~=~\sum\limits^k_{j=-\infty}z^jD~ \to D(k)_i~=~\sum\limits^k_{j=-\infty}z^jD_i~,\cr
&\phi_F(k)~=~\pmatrix  0 & -zh_F & 0 & \dots \cr 
               0 & 0 & -zh_F & \dots \cr 
               0 & 0 & 0 & \dots \cr
               \vdots & \vdots & \vdots & \ddots \endpmatrix~:\cr
&\hskip100pt
D(k)_i~=~\sum\limits^k_{j=-\infty}z^jD~ \to E(k)_i~=~\sum\limits^k_{j=-\infty}z^jE_i~.
\endaligned$$
As in 1.8 $\phi_D(k)$ is an automorphism, with inverse
$$\aligned
&\phi_D(k)^{-1}~=~\pmatrix 1 & zh_D & z^2(h_D)^2 & \dots \cr 
               0 & 1 & zh_D & \dots \cr 
               0 & 0 & 1 & \dots \cr
               \vdots & \vdots & \vdots & \ddots \endpmatrix~:\cr
&\hskip100pt
D(k)_i~=~\sum\limits^k_{j=-\infty}z^jD_i \to D(k)_i~=~\sum\limits^k_{j=-\infty}z^jD_i~,
\endaligned$$
and the chain complex $\overline{F}(k)$ defined by
$$\aligned
d_{\overline{F}(k)}~=~d_F-\phi_F(k)\phi_D(k)^{-1}c~&=~
\pmatrix d_F & zh_Fc & z^2h_Fh_Dc & \dots \cr 
               0 & d_F & zh_Fc & \dots \cr 
               0 & 0 & d_F & \dots \cr
               \vdots & \vdots & \vdots & \ddots \endpmatrix~:\cr
&\overline{F}(k)_i~=~\sum\limits^k_{j=-\infty}z^jF_i~ \to 
\overline{F}(k)_i~=~\sum\limits^k_{j=-\infty}z^jF_i
\endaligned$$
is such that the inclusions $\overline{F}(k)_i \to E(k)_i$ induce isomorphisms
$$\overline{F}(k)~\cong~F(k)~~,~~\varprojlim_k\overline{F}(k)~\cong~
\varprojlim_kF(k)~.$$
The identification 
$\displaystyle{\varprojlim_k}D(k)=D_{\alpha }((z))$ 
is immediate from the identifications
$$D(k)~=~\sum\limits^k_{j=-\infty}z^jD~ \to D(k-1)~=~
\sum\limits^{k-1}_{j=-\infty}z^jD~;~\sum\limits^k_{j=-\infty}z^jx_j \mapsto 
\sum\limits^{k-1}_{j=-\infty}z^jx_j~.$$
Similarly for 
$$\displaystyle{\varprojlim_k}E(k)~=~E_{\alpha }((z))~~,~~
\displaystyle{\varprojlim_k}\overline{F}(k)~=~\overline{F}~.$$
The short exact sequence of inverse systems
$$0 \to D(k) \raise4pt\hbox{$\phi\atop\to$} E(k) \to F(k) \to 0$$
determines an exact sequence of the inverse and derived limits
$$0 \to \varprojlim_kD(k) \raise4pt\hbox{$g-zh\atop\to$}  \varprojlim_kE(k) \to \varprojlim_kF(k) \to 
{\varprojlim_k}^1 D(k) \to \dots~.$$
Since the structure maps $D(k) \to D(k-1)$ are onto the derived limit is
$${\varprojlim_k}^1 D(k)~=~0$$
and the inverse limits actually fit into a short exact sequence
$$0 \to \varprojlim_kD(k) \raise4pt\hbox{$\overline{\phi}\atop\to$}  \varprojlim_kE(k) \to \varprojlim_kF(k) \to 0$$
as required, with $\overline\phi=g-zh$. Alternatively, identify
$\hbox{\rm coker}(\overline{\phi})=\overline{F}$ by a direct 
application of 1.8 (i).\qed
\enddemo

The following result on algebraic handle exchanges will be required in \S6.
$$\beginpicture
\setcoordinatesystem units <8pt,10pt>  
\putrule from -17 0 to 17 0
\putrule from -17 -5 to 17 -5
\putrule from -17 0 to -15 -5
\putrule from -15 0 to -15 -5
\putrule from -5 0 to -5 -5
\putrule from 5 0 to 5 -5
\putrule from 15 0 to 15 -5
\put {$E^+$} at -10 -2.7
\put {$E^-$} at 0 -2.7
\put {$\alpha E^+$} at 10 -2.7 
\put {$D$} at -15 -6
\put {$D'$} at -5 -6
\put {$\alpha D$} at 5 -6 
\put {$\alpha D'$} at 15 -6
\put {$\dots$} at  20 0
\put {$\dots$} at  20 -5
\put {$\dots$} at  -20 0
\put {$\dots$} at  -20 -5
\endpicture$$

\bigskip
\noindent{\bf Proposition 5.4} {\it Let $D,D',E^+,E^-$ be f.g. 
free $A$-module chain complexes, and let 
$$g^+~:~D \to E^+~,~g^-~:~D' \to E^-~,~
h^+~:~D' \to E^+~,~h^-~:~\alpha D  \to E^-$$
be $A$-module chain maps such that
$$\aligned
&g^+~=~\pmatrix 1 \cr 0 \endpmatrix~:~D_i \to E^+_i ~=~D_i \oplus F^+_i~,\cr
&g^-~=~\pmatrix 1 \cr 0 \endpmatrix~:~D'_i \to E^-_i ~=~D'_i \oplus F^-_i~,\cr
&h^+~=~\pmatrix h^+_D \cr h^+_F \endpmatrix~:~
D'_i \to E^+_i ~=~D_i \oplus F^+_i~,\cr
&h^-~=~\pmatrix h^-_D \cr h^-_F \endpmatrix~:~
\alpha D_i \to E^-_i ~=~D'_i \oplus F^-_i~,\cr
&d_{E^+}~=~\pmatrix d_D & c^+ \cr 0 & d_{F^+} \endpmatrix~:~
E^+_i~=~D_i \oplus F^+_i \to E_{i-1}~=~D_{i-1} \oplus F^+_{i-1}~,\cr
&d_{E^-}~=~\pmatrix d_{D'} & c^- \cr 0 & d_{F^-} \endpmatrix~:~
E^-_i~=~D'_i \oplus F^-_i \to E^-_{i-1}~=~D'_{i-1} \oplus F^-_{i-1}~.
\endaligned$$
Given bases for $D,D',F^+,F^-$ 
there are now defined two collections of data as in 5.3
$$(g:D \to E,h:\alpha D \to E,F)~~,~~(g':D' \to E',h':\alpha D' \to E',F')$$
with 
$$\aligned
&d_E~=~\pmatrix d_D & c \cr 0 & d_F \endpmatrix~:~
E_i~=~D_i \oplus F_i \to E_{i-1}~=~D_{i-1} \oplus F_{i-1}~,\cr
&d_F~=~\pmatrix d_{F^+} & h^+_Fc^- \cr 0 & d_{F^-} \endpmatrix~:~
F_i~=~F^+_i \oplus F^-_i \to F_{i-1}~=~F^+_{i-1} \oplus F^-_{i-1}~,\cr
&c~=~\pmatrix c^+ & h^+_Fc^- \endpmatrix~:~F_i~=~F^+_i \oplus F^-_i \to
D_{i-1}~,\cr
&g~=~\pmatrix 1 \cr  0 \endpmatrix~:~D_i \to E_i~=~D_i \oplus F_i~,\cr
&h_D~=~h^+_Dh^-_D ~:~\alpha D_i \to D_i~,\cr
&h_F~=~\pmatrix h^+_Fh^-_D \cr h^-_F \endpmatrix ~:~
\alpha D_i \to F_i~=~F^+_i \oplus F^-_i
\endaligned$$
and
$$\aligned
&d_{E'}~=~\pmatrix d_{D'} & c' \cr 0 & d_{F'} \endpmatrix~:~
E'_i~=~D'_i \oplus F'_i \to E'_{i-1}~=~D'_{i-1} \oplus F'_{i-1}~,\cr
&d_{F'}~=~\pmatrix d_{F^+} & 0 \cr 
h^-_Fc^+ & d_{F^-} \endpmatrix~:~
F'_i~=~F^+_i \oplus F^-_i \to F'_{i-1}~=~F^+_{i-1} \oplus F^-_{i-1}~,\cr
&c'~=~\pmatrix h^-_Dc^+ & c^- \endpmatrix~:~
F'_i~=~F^+_i \oplus F^-_i \to D'_{i-1}~,\cr
&g'~=~\pmatrix 1 \cr  0 \endpmatrix~:~D'_i \to E'_i~=~D'_i \oplus F'_i~,\cr
&h'_D~=~h^-_Dh^+_D ~:~\alpha D'_i \to D'_i~,\cr
&h'_F~=~\pmatrix h^+_F \cr h^-_Fh^+_D \endpmatrix ~:~
\alpha D'_i \to F'_i~=~F^+_i \oplus F^-_i~.
\endaligned$$
The cokernels of the corresponding embeddings of the based
f.g. free $\Sigma^{-1}A_{\alpha }[z,z^{-1}]$-module chain complexes
$$\aligned
\phi~=~g-zh~:~\Sigma^{-1}D_{\alpha }[z,z^{-1}] \to \Sigma^{-1}E_{\alpha }[z,z^{-1}]~,\cr
\phi'~=~g'-zh'~:~\Sigma^{-1}D'_{\alpha }[z,z^{-1}] \to \Sigma^{-1}E'_{\alpha }[z,z^{-1}]
\endaligned$$
are related by an isomorphism of the based f.g. 
free $\Sigma^{-1}A_{\alpha }[z,z^{-1}]$-module chain complexes
$$I~:~\hbox{\rm coker}(\phi)~\cong~\hbox{\rm coker}(\phi')$$
which sends basis elements to $z^\delta$(basis elements), with $\delta=0$ or $1$.}
\demo{Proof} Use 5.3 to identify
$$\hbox{\rm coker}(\phi)~=~\widehat{F}~~,~~\hbox{\rm coker}(\phi')~=~
\widehat{F'}$$
with $\widehat{F},\widehat{F}'$ the based f.g. free 
$\Sigma^{-1}A_{\alpha }[z,z^{-1}]$-module chain complexes defined by
$$\aligned
&d_{\widehat{F}}~=~d_F + zh_F(1-zh_D)^{-1}c~:\cr
&\hskip50pt 
\widehat{F}_i~=~\Sigma^{-1}(F_i)_{\alpha }[z,z^{-1}] \to
\widehat{F}_{i-1}~=~\Sigma^{-1}(F_{i-1})_{\alpha }[z,z^{-1}]~,\cr
&d_{\widehat{F}'}~=~d_{F'} + zh'_F(1-zh'_D)^{-1}c'~:\cr
&\hskip50pt
\widehat{F}'_i~=~\Sigma^{-1}(F'_i)_{\alpha }[z,z^{-1}] \to
\widehat{F}'_{i-1}~=~\Sigma^{-1}(F'_{i-1})_{\alpha }[z,z^{-1}]~.
\endaligned$$
Define an isomorphism of based f.g. free 
$\Sigma^{-1}A_{\alpha }[z,z^{-1}]$-module chain complexes 
$$I~:~\widehat{F}~\cong~\widehat{F}'$$
by 
$$\aligned
&I~:~\widehat{F}_i~=~\Sigma^{-1}(F^+_i \oplus F^-_i)_{\alpha }[z,z^{-1}]
\to \widehat{F}'_i~=~\Sigma^{-1}(F^+_i \oplus F^-_i)_{\alpha }[z,z^{-1}]~;\cr
&\hskip125pt  (x^+,x^-) \to (x^+,zx^-)~.
\endaligned\eqno{\lower12pt\hbox{\qed}}$$
\enddemo

\subheading{6. The Novikov complexes $C^{Nov}(M,f,v)$, $C^{Paj}(M,f,v)$, $C^{FR}(M,f,v,h)$} 

This section starts with a review of the geometric constructions of the
Novikov complex $C^{Nov}(M,f,v)$ and the Pajitnov complex $C^{Paj}(M,f,v)$
of a circle-valued Morse function $f:M\to S^1$ with respect to $v \in \GT(f)$.
Then the proper real-valued Morse function
$$\overline{f}~:~\overline{M}~=~f^*\R \to \R$$
is used to identify the Novikov complex with the inverse limit of the
proper Morse-Smale complexes 
$$C^{Nov}(M,f,v)~=~\varprojlim_k C^{MS}(M(k),f(k),v(k))$$
of the proper real-valued Morse functions
$$f(k)~=~\overline{f}\vert~:~M(k)~=~
\overline{f}^{\,-1}[-k,\infty) \to [-k,\infty)~~(k \geq 0)$$
with $v(k)=\overline{v}\vert$.  This is followed by a review of the
algebraic construction of the chain complex $C^{FR}(M,f,v,h)$
of Farber and Ranicki \cite{FR}. Finally, all this is put together to prove the
Cokernel, Invariance and Isomorphism Theorems already stated in
the Introduction. In particular, for $v\in \GT(f)$ with a gradient-like 
chain approximation $h^{gra}$ there exist basis-preserving isomorphisms
$$\aligned
&C^{Nov}(M,f,v)~\cong~C^{FR}(M,f,v,h^{gra};\widehat{\Z[\pi_1(M)]})~,\cr
&C^{Paj}(M,f,v)~\cong~C^{FR}(M,f,v,h^{gra})~.
\endaligned$$
\indent The infinite cyclic cover of $M$ determined by $f:M \to S^1$ is
$$\overline{M}~=~f^*\R~=~\{(x,y) \in M \times \R\,\vert\,f(x)=[y] \in S^1\}~,$$
with
$$\overline{f}~:~\overline{M} \to \R~;~(x,y) \mapsto y~.$$
The generating covering translation 
$$z~:~\overline{M} \to \overline{M}~;~(x,y) \mapsto (x,y-1)~.$$
is parallel to the downward $v$-gradient flow (and so acts from right to left,
or rather from top to bottom). Assume that $M$ and $\overline{M}=f^*\R$ are connected, 
so that 
$$\pi_1(M)~=~\pi_1(\overline{M}) \times_{\alpha } \Z$$
with 
$$\alpha ~=~z_*~:~\pi_1(\overline{M}) \to \pi_1(\overline{M})$$ 
the monodromy automorphism. The group ring of $\pi_1(M)$ is the
$\alpha $-twisted Laurent polynomial extension of $\Z[\pi_1(\overline{M})]$
$$\Z[\pi_1(M)]~=~\Z[\pi_1(\overline{M})]_{\alpha }[z,z^{-1}]~.$$
Write the Novikov ring of $\Z[\pi_1(M)]$ as
$$\widehat{\Z[\pi_1(M)]}~=~\Z[\pi_1(\overline{M})]_{\alpha }((z))~=~
\Z[\pi_1(\overline{M})]_{\alpha }[[z]][z^{-1}]~.$$

\noindent{\bf Definition 6.1} (\cite{N}, \cite{P1})
The {\it Novikov complex} $C^{Nov}(M,f,v)$ of a Morse function $f:M \to S^1$
with respect to $v \in \GT(f)$
is the based f.\,g. free $\widehat{\Z[\pi_1(M)]}$-module chain complex with
\roster
\item"{\rm (i)}" 
$\hbox{\rm rank}_{\widehat{\Z[\pi_1(M)]}}C^{Nov}_i(M,f,v)=c_i(f),$
with one basis element $\widetilde{p}$ for each critical
point $p \in M$ of index $i$, corresponding to a choice
of lift $\widetilde{p} \in \widetilde{M}$
\item"{\rm (ii)}" the boundary $\widehat{\Z[\pi_1(M)]}$-module morphisms are given by
$$\aligned
d~:~C^{Nov}_i(M,f,v)~=~\widehat{\Z[\pi_1(M)]}^{c_i(f)}&\to 
C^{Nov}_{i-1}(M,f,v)~=~\widehat{\Z[\pi_1(M)]}^{c_{i-1}(f)}~;\cr
&\widetilde{p} \mapsto \sum\limits_{\widetilde{q}}
\sum\limits_{u \in \pi_1(M)} n(\widetilde{p},u\widetilde{q})u\widetilde{q}
\endaligned$$
with $n(\widetilde{p},u\widetilde{q}) \in \Z$ the algebraic number of 
$\widetilde{v}$-gradient flow lines in $\widetilde{M}$ from
$\widetilde{p}$ to $u\widetilde{q}$.\qed
\endroster

\medskip

Let $\Sigma^{-1}\Z[\pi_1(M)]$ be the localization of $\Z[\pi_1(M)]$
defined in 5.2.

\noindent{\bf Definition 6.2} (Pajitnov \cite{P2,P3,P4})\newline
For a Morse function $f:M \to S^1$ and $v\in\GCCT(f)$ the Novikov
complex is of the form
$$C^{Nov}(M,f,v)~=~\widehat{\Z[\pi_1(M)]}\otimes_{\Sigma^{-1}\Z[\pi_1(M)]}C^{Paj}(M,f,v)$$
with the {\it Pajitnov complex} $C^{Paj}(M,f,v)$ a based f.\,g. free 
$\Sigma^{-1}\Z[\pi_1(M)]$-module chain such that
$$\hbox{\rm rank}_{\Sigma^{-1}\Z[\pi_1(M)]}C^{Paj}_i(M,f,v)~=~c_i(f)~.\qed$$

From now on it will be assumed that $0 \in S^1$ is a regular
value of $f$, with inverse image a codimension 1 framed submanifold
$$N^{m-1}~=~f^{-1}(0) \subset M^m~.$$
Thus $\overline{f}:\overline{M} \to \R$ is transverse regular 
at $\Z \subset \R$, and cutting $M$ along $N$ gives a cobordism
$$(M_N;N,z^{-1}N)~=~\overline{f}^{\,-1}(I;\{0\},\{1\})$$ 
which is a fundamental domain for $\overline{M}$
$$\overline{M}~=~\bigcup\limits^{\infty}_{j=-\infty}z^jM_N$$
with a Morse function
$$f_N~=~\overline{f}\vert~:~(M_N;N,z^{-1}N) \to (I;\{0\},\{1\})~.
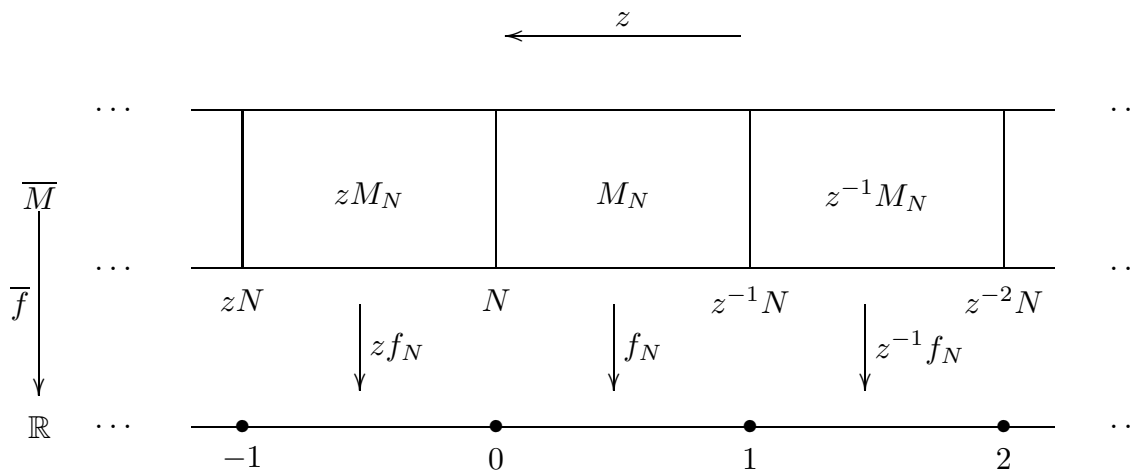
$$
$$\beginpicture
\setcoordinatesystem units <8pt,10pt>  
\putrule from -17 0 to 17 0
\putrule from -17 -5 to 17 -5
\putrule from -15 0 to -15 -5
\putrule from -5 0 to -5 -5
\putrule from 5 0 to 5 -5
\putrule from 15 0 to 15 -5
\putrule from -17 -10 to 17 -10
\put {$zM_N$} at -10 -2.7
\put {$M_N$} at 0 -2.7
\put {$z^{-1}M_N$} at 10 -2.7 
\put {$zN$} at -15 -6
\put {$N$} at -5 -6
\put {$z^{-1}N$} at 5 -6 
\put {$z^{-2}N$} at 15 -6
\put {$-1$} at -15 -11
\put {0} at -5 -11
\put {1} at 5 -11
\put {2} at 15 -11
\put {$\bullet$} at -15 -10
\put {$\bullet$} at -5 -10
\put {$\bullet$} at 5 -10
\put {$\bullet$} at 15 -10
\put {$\overline{M}$} at -23 -2.7
\put {$\R$} at -23 -10
\put {$\xymatrix@R+34pt {\ar[d]_{\displaystyle \overline{f}} & \\ &}$} at -21.7 -6.1
\put {$\xymatrix@R+3pt {\ar[d]^{\displaystyle {zf_N}} & \\ &}$} at -8.5 -7.5
\put {$\xymatrix@R+3pt {\ar[d]^{\displaystyle {f_N}} & \\ &}$} at 1.5 -7.5
\put {$\xymatrix@R+3pt {\ar[d]^{\displaystyle {z^{-1}f_N}} & \\ &}$} at 11.5 -7.5
\put {$\xymatrix@C+50pt { & \ar[l]_{\displaystyle z}} $} at 0 2.7
\put {$\dots$} at  20 0
\put {$\dots$} at  20 -5
\put {$\dots$} at  20 -10
\put {$\dots$} at  -20 0
\put {$\dots$} at  -20 -5
\put {$\dots$} at  -20 -10
\endpicture$$
\vskip2pt

\noindent{\bf Proposition 6.3} {\it The Novikov complex of a Morse function
$f:M \to S^1$ with respect to any $v \in \GT(f)$ is the inverse limit
$$C^{Nov}(M,f,v)~=~\varprojlim_k C^{MS}(M(k),f(k),v(k))$$
of the Morse-Smale complexes of the proper real-valued Morse functions
$$f(k)~=~\overline{f}\vert~:~M(k)~=~
\overline{f}^{\,-1}[-k,\infty)~=~\bigcup^k_{j=-\infty}z^jM_N
 \to [-k,\infty)$$
on the non-compact manifolds with boundary
$$(M(k),\partial M(k))~=~
\overline{f}^{\,-1}([-k,\infty),\{-k\})~=~(\bigcup^k_{j=-\infty}z^jM_N,z^kN)$$
with respect to the projections}
$$\aligned
&C^{MS}(M(k+1),f(k+1),v(k+1))~=~C(\bigcup^{k+1}_{j=-\infty}z^jM_N,z^{k+1}N)\\
&\to C(\bigcup^{k+1}_{j=-\infty}z^jM_N,z^{k+1}M_N)~=~
C(\bigcup^k_{j=-\infty}z^jM_N,z^kN)~=~C^{MS}(M(k),f(k),v(k))~.
\endaligned$$
\demo{Proof} Identify 
$$C^{Nov}(M,f,v)~=~\varprojlim_k C^{MS}(M(k),f(k),v(k))$$
using the expression of the Novikov ring as the inverse limit
$$\widehat{\Z[\pi_1(M)]}~=~\varprojlim_k \Z[\pi_1(\overline{M})]_{\alpha }(k)$$
of the $\Z[\pi_1(\overline{M})]$-modules
$$\Z[\pi_1(\overline{M})]_{\alpha }(k)~=~\sum\limits^k_{j=-\infty}z^j \Z[\pi_1(\overline{M})]
\subset \Z[\pi_1(M)]~=~\Z[\pi_1(\overline{M})]_{\alpha }[z,z^{-1}]$$
with respect to the natural projections
$$\Z[\pi_1(\overline{M})]_{\alpha }(k+1) \to
\Z[\pi_1(\overline{M})]_{\alpha }(k)~;~
\sum\limits^k_{j=-\infty}a_jz^j \mapsto \sum\limits^{k-1}_{j=-\infty}a_jz^j~.\qed$$
\enddemo

As before, let $f:M \to S^1$ have $c_i(f)$ critical points of index $i$.
The real-valued Morse function
$$f_N~=~\overline{f}\vert~:~(M_N;N,z^{-1}N) \to (I;\{0\},\{1\})$$
has $c_i(f_N)=c_i(f)$ critical points of index $i$, and as in \S2 
there is a handle decomposition
$$M_N~=~N \times I \cup \bigcup\limits^m_{i=0}\bigcup\limits_{c_i(f)}h^i~.$$
Let $N$ have a $CW$ structure with 
$c_i(N)$ $i$-cells $e^i \subset N$ and let $M_N$ have the corresponding $CW$ structure
with $c_i(N)$ $i$-cells of type $e^i \times I \subset M_N$ and $c_i(f)$
$i$-cells of type $h^i \subset M_N$. Let $\widetilde{M}$ be the universal
cover of $M$, and let $\widetilde{M}_N,\widetilde{N}$ be the corresponding
covers of $M_N,N$. The inclusion $g:N \to M_N$ induces an inclusion of
the cellular $\Z[\pi_1(\overline{M})]$-module chain complexes 
$g: C(\widetilde{N}) \to C(\widetilde{M}_N)$. The inclusion
$h:z^{-1}N \to M_N$ 
induces a $\Z[\pi_1(\overline{M})]$-module chain map
$h:\alpha C(\widetilde{N}) \to C(\widetilde{M}_N)$.
$$\beginpicture
\setcoordinatesystem units <6pt,10pt>  
\putrule from 0 0 to 20 0
\putrule from 0 5 to 20 5
\putrule from 0 0 to 0 5
\putrule from 20 0 to 20 5
\put {$M_N$} at 10  2.7
\put {$N$} at -1.5 2.7
\put {$z^{-1}N$} at 23 2.7
\put {$\xymatrix@C+10pt {\ar[r]^{\displaystyle{g}}& } $ } at 4.2 3
\put {$\xymatrix@C+10pt {& \ar[l]_{\displaystyle{h}}}$ } at 16 3
\endpicture$$
\vskip2pt

\noindent Now $M=M_N/(N=z^{-1}N)$ has a $CW$ structure with 
$$c_i(M)~=~c_i(f)+c_i(N)+c_{i-1}(N)$$ 
$i$-cells.

The following terminology will be used in dealing with the
cellular chain complexes associated to the fundamental domain
$(M_N;N,z^{-1}N)$ of $\overline{M}$.

\noindent{\bf Terminology 6.4} 
As in 2.3 write
$$\aligned
&D~=~C(\widetilde{N})~~,~~E~=~C(\widetilde{M}_N)~~,~~
F~=~C^{MS}(M_N,f_N,v_N)~=~C(\widetilde{M}_N,\widetilde{N}\times I)~,\cr
&d_E~=~\pmatrix d_D & c \cr 0 & d_F \endpmatrix~:~E_i~=~D_i \oplus F_i \to
E_{i-1}~=~D_{i-1} \oplus F_{i-1}~.
\endaligned$$
The inclusions 
$$g:N \to M~~,~~h:z^{-1}N \to M$$
induce chain maps
$$\aligned
&g~=~\pmatrix 1 \cr 0 \endpmatrix~:~
D_i \to E_i~=~D_i\oplus F_i ~,\cr
&h~=~\pmatrix h_D \cr h_F\endpmatrix~:~
\alpha D_i \to E_i~=~D_i\oplus F_i~.\qed
\endaligned
$$

\noindent{\bf Definition 6.5} (\cite{FR})
Given a Morse function $f:M \to S^1$ with regular value $0 \in S^1$,
$v \in \GT(f)$, a choice of $CW$ structure for $N=f^{-1}(0) \subset M$, 
and a choice of chain approximation
$h:\alpha C(\widetilde{N}) \to C(\widetilde{M}_N)$ 
$$C^{FR}(M,f,v,h)~=~\widehat{F}$$ 
be the based f.\,g. free $\Sigma^{-1}\Z[\pi_1(M)]$-module chain complex 
given by
$$\aligned
&\widehat{F}_i~=~\Sigma^{-1}C_i(\widetilde{M}_N,\widetilde{N})_{\alpha }[z,z^{-1}]~=~
\Sigma^{-1}\Z[\pi_1(M)]^{c_i(f)}~,\\
&d_{\widehat{F}}~=~d_{C(\widetilde{M}_N,\widetilde{N})}+zh_F(1-zh_D)^{-1}c~:~
\widehat{F}_i \to \widehat{F}_{i-1}~.\qed
\endaligned$$

\proclaim{\bf Cokernel Theorem 6.6} 
{\rm (i)} The inclusions
$\widehat{F}_i \to \Sigma^{-1}C_i(\widetilde{M}_N)_{\alpha }[z,z^{-1}]$
induce a basis-preserving {\rm (}and a fortiori simple{\rm )}
isomorphism of based f.\,g. free $\Sigma^{-1}\Z[\pi_1(M)]$-module chain complexes
$$C^{FR}(M,f,v,h) ~\cong~\hbox{\rm coker}
\big(g-zh:\Sigma^{-1}C(\widetilde{N})_{\alpha }[z,z^{-1}] \to 
\Sigma^{-1}C(\widetilde{M}_N)_{\alpha }[z,z^{-1}]\big)~.$$
{\rm (ii)} The natural projection
$$\aligned
&p~:~C(\widetilde{M};\Sigma^{-1}\Z[\pi_1(M)])~=~
{\Cal C}\big(g-zh:\Sigma^{-1}C(\widetilde{N})_{\alpha }[z,z^{-1}] \to 
\Sigma^{-1}C(\widetilde{M}_N)_{\alpha }[z,z^{-1}]\big)\cr
&\hskip75mm \to C^{FR}(M,f,v,h)
\endaligned$$
is a chain equivalence of based f.\,g. free 
$\Sigma^{-1}\Z[\pi_1(M)]$-module chain complexes with torsion
$$\aligned
&\tau(p)~=~-\sum\limits^{\infty}_{i=0}(-)^i\tau\big(1-zh_D:
\Sigma^{-1}C_i(\widetilde{N})_{\alpha }[z,z^{-1}] \to 
\Sigma^{-1}C_i(\widetilde{N})_{\alpha }[z,z^{-1}]\big)\cr
&\hskip50mm  \in K_1(\Sigma^{-1}\Z[\pi_1(M)])~.
\endaligned$$
\endproclaim
\demo{Proof} This is a direct application of 5.3 (i). \hfill\qed
\enddemo

\noindent{\bf Invariance Theorem 6.7} 
{\it Let $f:M \to S^1$ be a Morse function, and let $v \in \GT(f)$.
For any regular values $0,0' \in S^1$, 
$CW$ structures on $N=f^{-1}(0)$, $N'=f^{-1}(0')$ and chain approximations
$$h~:~\alpha C(\widetilde{N}) \to C(\widetilde{M}_N)~,~
h'~:~\alpha C(\widetilde{N}') \to C(\widetilde{M}_{N'})$$
there is defined a simple isomorphism of 
based f.\,g. free $\Sigma^{-1}\Z[\pi_1(M)]$-module chain complexes
$$C^{FR}(M,f,v,h)~\cong~C^{FR}(M,f,v,h')~.$$}
\demo{Proof} The case $0=0'$, $N=N'$ is given by 5.3 (ii). 
So assume $0 \not = 0' \in S^1$, and let $I^+,I^- \subset S^1$ 
be the two arcs joining 0 and $0'$. The restrictions of $(f,v)$
$$\aligned
&(f^+,v^+)~=~(f,v)\vert~:~(M^+;N,N')~=~f^{-1}(I^+;\{0\},\{0'\}) \to (I^+;\{0\},\{0'\})~,\cr
&(f^-,v^-)~=~(f,v)\vert~:~(M^-;N',N)~=~f^{-1}(I^-;\{0'\},\{0\}) \to (I^-;\{0'\},\{0\})
\endaligned$$
are Morse functions with gradient-like vector fields such that
$$\aligned
&(f_N,v_N)~=~(f^+,v^+)\cup (f^-,v^-)~:\cr
&(M_N;N,z^{-1}N)~=~
(M^+;N,N') \cup (M^-;N',N) \to ([0,1];\{0\},\{1\})~,\cr
&(f_{N'},v_{N'})~=~(f^-,v^-)\cup (f^+,v^+)~:\cr
&(M_{N'};N',z^{-1}N')~=~
(M^-;N',N) \cup (M^+;N,N') \to ([0,1];\{0\},\{1\})~.
\endaligned$$
Use the handlebody structure on $(M^+;N,N')$ (resp.  $(M^-;N',N)$)
determined by $(f^+,v^+)$ (resp.  $(f^-,v^-)$) to extend the $CW$
structure on $N$ (resp.  $N'$) to a $CW$ structure on $M^+$ (resp. 
$M^-$).  The inclusions of $CW$ subcomplexes $g^+:N \to M^+$, 
$g^-:N' \to M^-$ induce inclusions of subcomplexes
$$g^+~:~D~=~C(\widetilde{N}) \to E^+~=~C(\widetilde{M}^+)~,~
g^-~:~D'~=~C(\widetilde{N}') \to E^+~=~C(\widetilde{M}^-)~.$$
Cellular approximations to the
inclusions $h^+:N' \to M^+$, $h^-:z^{-1}N\to M^-$ induce chain maps
$$h^+~:~D'~=~C(\widetilde{N}') \to E^+~=~C(\widetilde{M}^+)~,~
h^-~:~\alpha D~=~\alpha C(\widetilde{N}) \to E^-~=~C(\widetilde{M}^-)~.$$
A direct application of 5.4 gives a simple isomorphism of based f.g. free
$\Sigma^{-1}\Z[\pi_1(M)]$-module chain complexes
$$C^{FR}(M,f,v,h)~=~\hbox{\rm coker}(g-zh)~\cong~
C^{FR}(M,f,v,h')~=~\hbox{\rm coker}(g'-zh')~.\eqno{\qed}$$
\enddemo

\noindent{\bf Definition 6.8} 
A crossing chain approximation
$h:\alpha C(\widetilde{N}) \to C(\widetilde{M}_N)$ is 
{\it gradient-like} if for any critical points $\widetilde{p},\widetilde{q}
\in \widetilde{M}_N$ of $\widetilde{f}_N:\widetilde{M}_N \to \R$ with
index $i,i-1$ the algebraic number of $\widetilde{v}$-gradient 
flow lines in $\widetilde{M}$ joining $\widetilde{p}$ to $z^j \widetilde{q}$ is
$$n(\widetilde{p},z^j\widetilde{q})~=~
\cases
\hbox{\rm $(\widetilde{p},\widetilde{q})$-coefficient of $d_F:F_i \to F_{i-1}$}
&\hbox{\rm if $j=0$} \cr
\hbox{\rm $(\widetilde{p},\widetilde{q})$-coefficient of $h_F(h_D)^{j-1}c:F_i \to F_{i-1}$}
&\hbox{\rm if $j>0$} \cr
0&\hbox{\rm if $j<0$~.}
\endcases$$
A gradient-like chain approximation $h$ will be denoted by $h^{gra}$.\hfill\qed

\noindent{\bf Remark 6.9} 
(i) Pajitnov \cite{P4} proved that for any Morse function 
$f:M \to S^1$ and $v \in \GCCT(f)$ there exists a 
handlebody structure on $N$ with a gradient-like chain approximation
$h^{gra}:\alpha C(\widetilde{N}) \to C(\widetilde{M}_N)$
(cf. Remark 2.4).\newline
(ii) Cornea and Ranicki \cite{CR} prove that for any Morse function
$f:M \to S^1$ and $v \in \GT(f)$ there exists a $CW$ structure on $N=f^{-1}(0)$
with a gradient-like chain approximation $h^{gra}$ by showing that there exist Morse
functions $f':M \to S^1$, $g:N \to \R$ with $v' \in \GT(f')$, $w \in
\GT(g)$ such that
\roster
\item"(a)" $(f',v')$  agrees with $(f,v)$ outside a tubular neighbourhood of $N$
$$f^{-1}[-\epsilon,\epsilon]~=~N \times [-\epsilon,\epsilon] \subset M$$
for some small $\epsilon>0$.
\item"(b)" $(f',v')$ restricts to translates of $(g,w)$
$$\aligned
&(f',v')\vert~=~(g_+,w_+)~:~N \times \{\epsilon/2\} \to S^1 \backslash \{0\}~,\cr
&(f',v')\vert~=~(g_-,w_-)~:~N \times \{-\epsilon/2\} \to S^1 \backslash \{0\}
\endaligned$$
with 
$$\aligned
&\hbox{\rm Crit}_i(g_+)~=~\hbox{\rm Crit}_i(g) \times \{\epsilon/2\}~,\cr
&\hbox{\rm Crit}_i(g_-)~=~\hbox{\rm Crit}_i(g) \times \{-\epsilon/2\}~,\cr
&\hbox{\rm Crit}_i(f')~=~\hbox{\rm Crit}_{i-1}(g_+) 
\cup \hbox{\rm Crit}_i(g_-) \cup \hbox{\rm Crit}_i(f)~.
\endaligned$$
\item"(c)" The $\overline{v}$-gradient flow lines are in one-one
correspondence with the broken $\overline{v}'$-gradient flow lines i.e. 
joined up sequences of $\overline{v}'$-gradient flow lines which start
and terminate at critical points of $\overline{f}$.
\item"(d)" The Morse-Smale complex of $(\overline{f}',\overline{v}')$ is of the form
$$C^{MS}(\overline{M},\overline{f}',\overline{v}')~=~
((D_{i-1} \oplus D_i \oplus F_i)_{\alpha }[z,z^{-1}],
\pmatrix -d_D & 0 & 0 \cr
1-zh^{gra}_D & d_D & c \cr
-zh^{gra}_F &  0 & d_F \endpmatrix)$$
for a gradient-like chain approximation $h^{gra}$, with 
$$D~=~C^{MS}(N,g,w)~~,~~F~=~C^{MS}(M_N,f_N,v_N)~~\hbox{\rm etc.}$$
\item"(e)" The cellular chain complex of the universal cover $\widetilde{M}_N$
(or rather the cover of $M_N$ induced from the universal cover $\widetilde{M}$
of $\overline{M}$) is
$$E~=~C(\widetilde{M}_N)~=~(D_i \oplus F_i,\pmatrix d_D & c \cr
0 & d_F \endpmatrix)~.$$
\endroster
(This is the circle-valued analogue of Remark 2.7). Thus the
algebraic mapping cone of the $\Z[\pi_1(M)]$-module chain map
$$\phi~=~\pmatrix 1-zh_D^{gra} \cr -zh^{gra}_F \endpmatrix~:~
D_{\alpha}[z,z^{-1}] \to E_{\alpha}[z,z^{-1}]$$
is the Morse-Smale complex of 
$(\overline{f}',\overline{v}'):\overline{M} \to \R$
$${\Cal C}(\phi)~=~C^{MS}(\overline{M},\overline{f}',\overline{v}')$$
and the cokernel of the induced $\Sigma^{-1}\Z[\pi_1(M)]$-module chain map is
$$\hbox{\rm coker}(\Sigma^{-1}\phi:\Sigma^{-1}D_{\alpha}[z,z^{-1}] \to
\Sigma^{-1}E_{\alpha}[z,z^{-1}])~=~C^{FR}(M,f,v,h^{gra})~.$$
The kernel of the projection
$$K~=~\hbox{\rm ker}(p:\Sigma^{-1}{\Cal C}(\phi) \to \hbox{\rm coker}(\Sigma^{-1}\phi))$$
is an algebraic model for the closed orbits of the $v$-gradient flow.
As in 1.8 (i) (c) identify
$$\aligned
&d_K~=~\pmatrix d_D & 0 \cr (-1)^i(1-zh^{gra}_D) & (1-zh^{gra}_D)^{-1}d_D(1-zh^{gra}_D) \endpmatrix~:\cr
&K_i~=~\Sigma^{-1}(D_{i-1}\oplus D_i)_{\alpha}[z,z^{-1}] \to
K_{i-1}~=~\Sigma^{-1}(D_{i-2}\oplus D_{i-1})_{\alpha}[z,z^{-1}]~.
\endaligned$$
(iii) In applying 2.6 and 3.3 to $f:M \to S^1$, $v \in \GT(f)$ with 
a gradient-like chain approximation $h^{gra}$, the Morse-Smale complexes of
unions of copies of a Morse function 
$$f_N~:~(M_N;N;z^{-1}N) \to (I;\{0\},\{1\})$$
and $v_N \in \GT(f_N)$, the chain homotopy $b$ in 2.6
and the higher chain homotopies $b[j,j']$ in 3.3 are 0. In particular,
the Morse-Smale complex of the proper real-valued Morse function
$$\overline{f}^+~=~\bigcup^0_{k=-\infty}z^kf_N~=~\overline{f}\vert~:~
\overline{M}^+~=~\bigcup^0_{k=-\infty}z^kM_N~=~\overline{f}^{-1}[0,\infty)
\to [0,\infty)$$
is given as a based free $\Z[\pi_1(\overline{M})]$-module chain complex by
$$C^{MS}(\overline{M}^+,\overline{f}^+,\overline{v}^+)~=~
(\sum\limits^0_{k=-\infty}z^kF_i,d_F+\sum\limits^{\infty}_{j=1}z^jh^{gra}_F(h^{gra}_D)^{j-1}c)~.$$
The coefficients of $h^{gra}_F(h^{gra}_D)^{j-1}c$ count the 
$\widetilde{v}$-gradient flow
lines which start at an index $i$ critical point of $\widetilde{M}$,
cross $j$ translates of $\widetilde{N} \subset \widetilde{M}$, and
terminate at an index $i-1$ critical point of $\widetilde{M}$.\hfill\qed

\proclaim{Isomorphism Theorem 6.10} 
Given a Morse function $f:M \to S^1$ with regular value $0 \in S^1$,
$v \in \GT(f)$, a choice of $CW$ structure for $N=f^{-1}(0) \subset M$, 
and a choice of chain approximation
$h:\alpha C(\widetilde{N}) \to C(\widetilde{M}_N)$ 
there is defined a simple isomorphism of based f.\,g. free 
$R$-module chain complexes with 
$R=\widehat{\Z[\pi_1(M)]}$ {\rm (}resp. $\Sigma^{-1}\Z[\pi_1(M)]${\rm )}
$$\aligned
&I_h~:~C^{Nov}(M,f,v)~\cong~C^{FR}(M,f,v,h;\widehat{\Z[\pi_1(M)]})\cr
&(\hbox{\it resp.}~I_h~:~C^{Paj}(M,f,v)~\cong~C^{FR}(M,f,v,h))~.
\endaligned$$
For a gradient-like chain approximation
$h^{gra}:\alpha C(\widetilde{N}) \to C(\widetilde{M}_N)$ 
the simple isomorphisms $I_{h^{gra}}$ are basis-preserving. 
\endproclaim
\demo{Proof} By 6.3 the Novikov complex is the inverse limit
$$C^{Nov}(M,f,v)~=~\varprojlim_kF(k)$$
of the inverse system 
$$F(k)~=~C^{MS}(M(k),f(k),v(k))$$
of the Morse-Smale $\Z[\pi_1(\overline{M})]$-module chain complexes
of the proper real-valued Morse functions
$$f(k)~=~\overline{f}\vert~:~(M(k),\partial M(k))~=~
(\bigcup^k_{j=-\infty}z^jM_N,z^kN) \to ([-k,\infty),\{-k\})~.$$
By 3.3 there is defined an isomorphism of inverse systems
$$I_h(k)~:~F(k) \to \hbox{\rm coker}\big(g-zh:D(k) \to E(k)\big)$$
with
$$D(k)~=~\sum\limits^k_{j=-\infty}z^jC(\widetilde{N})~,~
E(k)~=~\sum\limits^k_{j=-\infty}z^jC(\widetilde{M}_N)~.$$
The inverse limits are given by 5.3 (iii) and 6.3
$$\varprojlim_k D(k)~=~C(\widetilde{N})_{\alpha }((z))~,~
\varprojlim_k E(k)~=~C(\widetilde{M}_N)_{\alpha }((z))$$
with 
$$\aligned
\varprojlim_k \hbox{\rm coker}\big(g-zh:D(k) \to E(k)\big)~
&=~\hbox{\rm coker}\big(g-zh:\varprojlim_k D(k) \to \varprojlim_k E(k)\big)\cr
&=~\hbox{\rm coker}\big(g-zh:C(\widetilde{N})_{\alpha }((z))
\to C(\widetilde{M}_N)_{\alpha }((z))\big)\cr
&\cong~C^{FR}(M,f,v,h;\widehat{\Z[\pi_1(M)]})~.
\endaligned$$
Define $I_h$ to be the induced isomorphism of inverse limits
$$\aligned
&I_h~=~\varprojlim_k I_h(k)~:~C^{Nov}(M,f,v)~=~\varprojlim_kF(k)\cr
&\hskip25pt
\to \varprojlim_k\hbox{\rm coker}\big(g-zh:D(k)\to E(k)\big)~\cong~
C^{FR}(M,f,v,h;\widehat{\Z[\pi_1(M)]})~.
\endaligned$$
For a gradient-like chain approximation $h^{gra}$ 3.3 gives basis-preserving 
identifications
$$\aligned
&d_{F(k)}~=~d_{C(\widetilde{M}_N,\widetilde{N})}+
\sum\limits^{\infty}_{j=1}z^jh_F(h^{gra}_D)^{j-1}c~:\cr
&F(k)_i~=~\sum\limits^k_{j=-\infty}z^jC_i(\widetilde{M}_N,\widetilde{N})
\to
F(k)_{i-1}~=~\sum\limits_{j=-\infty}^kz^jC_{i-1}(\widetilde{M}_N,\widetilde{N})
\endaligned$$
(as in 6.9 (iii)). Passing to the inverse limit as $k \to \infty$ gives
$$\aligned
&d_{C^{Nov}(M,f,v)}~=~d_{C(\widetilde{M}_N,\widetilde{N})}+
\sum\limits^{\infty}_{j=1}z^jh^{gra}_F(h^{gra}_D)^{j-1}c~:\cr
&C_i^{Nov}(M,f,v)~=~C_i(\widetilde{M}_N,\widetilde{N})_{\alpha }((z)) \to
C_{i-1}^{Nov}(M,f,v)~=~C_{i-1}(\widetilde{M}_N,\widetilde{N})_{\alpha }((z))
\endaligned$$
so that 
$$I_{h^{gra}}~:~C^{Nov}(M,f,v)~\cong~C^{FR}(M,f,v,h^{gra};\widehat{\Z[\pi_1(M)]})$$
is a basis-preserving isomorphism, with zero torsion. There exists a
chain homotopy
$$h~ \simeq~ h^{gra}~:~\alpha C(\widetilde{N}) \to C(\widetilde{M}_N)$$
so that by 5.3 (ii) 
$$\tau(I_h)~=~\tau(I_{h^{gra}})~=~0 \in K_1(\widehat{\Z[\pi_1(M)]})~.$$
Similarly for $\Sigma^{-1}\Z[\pi_1(M)]$-coefficients.\qed
\enddemo

\noindent{\bf Remark 6.11} 
{\rm (i)} The formulae given by 6.7 and 6.10
$$\aligned
&\tau\big(p:C(\widetilde{M};\Sigma^{-1}\Z[\pi_1(M)]) \to C^{Paj}(M,f,v)\big)\cr
&\hbox{\hskip20pt}=~-\sum\limits^{\infty}_{i=0}(-)^i\tau\big(1-zh^{gra}_D:
\Sigma^{-1}C_i(\widetilde{N})_{\alpha }[z,z^{-1}] \to 
\Sigma^{-1}C_i(\widetilde{N})_{\alpha }[z,z^{-1}]\big) \\
&\hbox{\hskip75mm} \in K_1(\Sigma^{-1}\Z[\pi_1(M)])~,\\
&\tau\big(\widehat{p}:C(\widetilde{M};\widehat{\Z[\pi_1(M)]}) \to C^{Nov}(M,f,v)\big)\cr
&\hbox{\hskip20pt}=~-\sum\limits^{\infty}_{i=0}(-)^i\tau\big(1-zh^{gra}_D:
C_i(\widetilde{N};\widehat{\Z[\pi_1(M)]}) \to
C_i(\widetilde{N};\widehat{\Z[\pi_1(M)]})\big)\\
&\hbox{\hskip75mm} \in K_1(\widehat{\Z[\pi_1(M)]})
\endaligned$$
are generalizations of the formulae of Hutchings and Lee \cite{HL} and
Pajitnov \cite{P5},\cite{P6} counting the critical points of $f:M \to S^1$, 
the $\zeta$-function of the closed orbits of the
gradient flow (corresponding to $h^{gra}_D$)
and the Reidemeister torsion of $M$. Sch\"utz \cite{Sch1},\cite{Sch2} 
extended these formulae to the closed orbits of a generic gradient flow
of a closed 1-form using Hochschild homology and a chain equivalence of 
the type $\widehat{p}:C(\widetilde{M};\widehat{\Z[\pi_1(M)]}) \to C^{Nov}(M,f,v)$.
\newline
{\rm (ii)} See Chapters 10,14,15 of \cite{R}
for the splitting theorems for the torsion groups 
$K_1(\widehat{\Z[\pi_1(M)]})$, $K_1(\Sigma^{-1}\Z[\pi_1(M)])$ in
the case $\alpha =1$ (which  extend to the case of arbitrary $\alpha $,
Pajitnov and Ranicki \cite{PR})  
and for the expressions of the torsions $\tau(1-zh_D)$ 
in terms of (noncommutative) characteristic polynomials. 
For any ring $A$ the classical Bass-Heller-Swan splitting
$$K_1(A[z,z^{-1}])~=~K_1(A)\oplus K_0(A)\oplus \widetilde{\text{Nil}}_0(A)
\oplus \widetilde{\text{Nil}}_0(A)$$
generalizes to splittings
$$\aligned
&K_1(\Sigma^{-1}A[z,z^{-1}])~=~K_1(A)\oplus K_0(A)\oplus 
\widetilde{\text{Nil}}_0(A) \oplus V(A)~,\\
&K_1\big(A((z))\big)~=~K_1(A)\oplus K_0(A)\oplus \widetilde{\text{Nil}}_0(A)
\oplus \widehat{V}(A)
\endaligned$$
with $\widehat{V}(A)\subseteq K_1\big(A((z))\big)$ 
the image of the multiplicative group of noncommutative 
Witt vectors $\widehat{W}(A)=1+zA[[z]]$, 
and $V(A)\subseteq K_1(\Sigma^{-1}A[z,z^{-1}])$ the image of the
subgroup of the noncommutative rational Witt vectors $W(A) \subseteq \widehat{W}(A)$
generated by $1+zA[z]$. (In \cite{R} it was claimed that
the natural surjections $W(A)^{ab}\to V(A)$,
$\widehat{W}(A)^{ab}\to\widehat{V}(A)$ are isomorphisms, but an explicit
counterexample was constructed in \cite{PR}).
The torsions of the chain equivalences $p,\widehat{p}$ are such that
$$\tau(p) \in V(A) \subseteq K_1(\Sigma^{-1}A[z,\alpha ])~~,~~
\tau(\widehat{p}) \in \widehat{V}(A) \subseteq K_1\big(A((z))\big)$$
respectively, with $A=\Z[\pi_1(\overline{M})]$.\newline
(iii) The natural map $\Sigma^{-1}A[z,z^{-1}] \to A((z))$ is injective
for a commutative ring $A$, since in this case $\Sigma^{-1}A[z,z^{-1}]$
is just the localization of $A[z,z^{-1}]$ inverting all the elements
of type $1+az \in A[z,z^{-1}]$ ($a \in A)$.
Sheiham \cite{Sh} has constructed an example of a noncommutative ring $A$
such that $\Sigma^{-1}A[z,z^{-1}] \to A((z))$ is not injective.\newline
(iv) Farber \cite{F} extended the construction of Farber and Ranicki
\cite{FR} to obtain an algebraic Novikov complex for a closed 1-form,
but did not relate it to the geometric Novikov complex.
\qed

\enddocument